\documentclass[11pt,a4paper]{article}
\pdfoutput=1
\newcommand{\s}[1]{{\textsf{\textbf{#1}}}}

\usepackage{graphicx}
\usepackage{amsmath}
\usepackage{grffile}
\usepackage{mathrsfs}
\usepackage{amsthm}
\usepackage{amsfonts}
\usepackage{subfigure}

\newcommand\grad{{\bf \nabla}}

\newcommand\nvec{{\bf n}}
\newcommand\xvec{{\bf x}}

\newcommand\uvec{{\bf u}}

\newcommand\Vvec{{\bf V}}
\newcommand\Qvec{{\bf Q}}
\newcommand\mvec{{\bf m}}
\newcommand\eps{{\epsilon}}

\newcommand{\Rr}{{\mathbb R}}

\title{\huge\s{Nematic equilibria on a two-dimensional annulus: defects and energies}}
 \date{\today}
\author{{\Large Alexander H. Lewis$^1$, Peter D. Howell$^1$, Dirk G. A. L. Aarts$^2$,}\\{\Large and Apala Majumdar$^{3*}$}}

\begin{document}

\maketitle

\begin{abstract} We study planar nematic equilibria on a two-dimensional  annulus with strong and weak tangent anchoring, within the Oseen-Frank and Landau-de Gennes theories for nematic liquid crystals. We analyse the defect-free state in the Oseen-Frank framework and obtain analytic stability criteria in terms of the elastic anisotropy, annular aspect ratio and anchoring strength. We consider radial and azimuthal perturbations of the defect-free state separately, which yields a complete stability diagram for the defect-free state. We construct nematic equilibria with an arbitrary number of defects on a two-dimensional annulus with strong tangent anchoring and compute their energies; these equilibria are generalizations of the diagonal and rotated states observed in a square. This gives novel insights into the correlation between preferred numbers of defects, their locations and the geometry. In the Landau-de Gennes framework, we adapt Mironescu's powerful stability result in the Ginzburg-Landau framework (P.~Mironescu, \textit{On the stability of radial solutions of the Ginzburg-Landau equation}, 1995) to compute quantitative criteria for the local stability of the defect-free state in terms of the temperature and geometry.

\end{abstract}

\footnotetext{\textit{$^{1}$Mathematical Institute, University of Oxford, Oxford, OX2 6GG, UK.}}
\footnotetext{\textit{$^{2}$Department of Chemistry, Physical and Theoretical Chemistry Laboratory, University of Oxford, Oxford, OX1 3QZ, UK.}}
\footnotetext{\it $^3$Department of Mathematical Sciences, University of Bath, Bath, BA2 7AY, UK}
\footnotetext{\it $^*$Corresponding author: A. M. (E-mail: \tt a.majumdar@bath.ac.uk)}
          
\section{Introduction}
\label{sec:intro}

\begin{figure}[ht]
\centering
\includegraphics[scale=1]{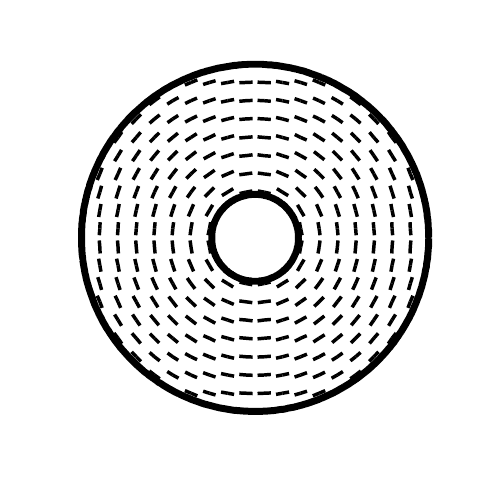}
\includegraphics[scale=1]{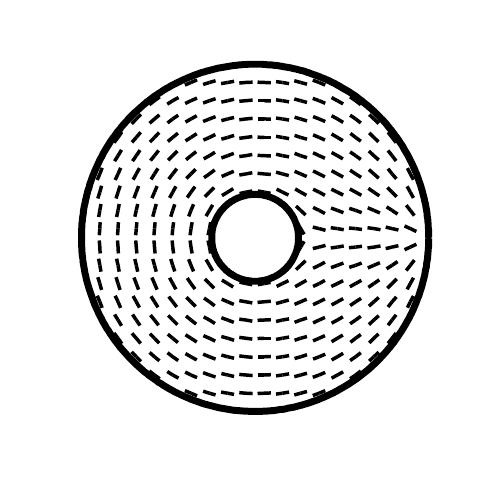}
\includegraphics[scale=1]{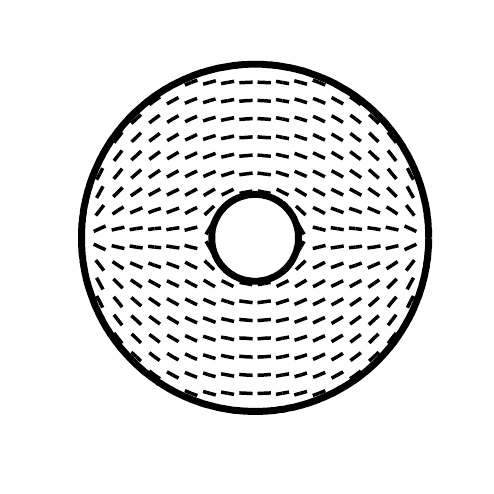}
\label{fig:Intro}
\caption{Nematic director fields identified experimentally in \cite{Oli_thesis, Jose_thesis}.}
\end{figure}

Nematic liquid crystals (LCs) are classic examples of partially ordered materials that combine the fluidity of liquids with the orientational order of crystalline solids \cite{dg,virga}.
Nematics have generated substantial scientific interest in recent years because of their unique optical, mechanical and rheological properties \cite{sluckin} and notably, nematics form the backbone of the multi-billion dollar liquid crystal display (LCD) industry. Defects are a key feature of nematic spatio-temporal patterns in confined geometries. A lot remains to be understood about the structure of defects and how they can be created, controlled and manipulated to yield desired properties. In this paper, we revisit the classical problem of a nematic sample in a two-dimensional (2D) annulus with strong or weak tangent boundary conditions separately and no external fields. Our work is motivated in part by recent experiments \cite{Oli_thesis, Jose_thesis}, on rod-like \textit{fd}-virus particles within shallow, annular microscopic chambers. Within these  chambers, multiple states are observed, including a radially invariant defect-free state and states with regularly arranged defects on the boundary, see Figure \ref{fig:Intro}.

The model problem of a defect-free state in annular wells has received a great deal of attention in the past, especially within the Oseen-Frank (OF) theory for liquid crystals. Here, we give a brief overview. The papers \cite{helein, pan, duffy1, duffy2, duffy3, palffy, virga2} are directly relevant to our work. In \cite{pan}, the author studies the stability and multiplicity of nematic radial equilibria on a 2D annulus with strong uniform anchoring. In \cite{duffy1, duffy2, palffy}, the authors approach the same problem with a more applications-oriented perspective motivated by the classical Freederickzs transition. They study confined nematic samples between two concentric cylinders with weak anchoring on one lateral surface and strong anchoring on another, subject to an external magnetic field. The authors primarily consider the stability of three characteristic configurations, referred to as `radial', `azimuthal' and `uniform' and obtain explicit estimates for the critical threshold field in terms of the OF elastic anisotropy, cylindrical aspect ratio and anchoring strength. The defect-free state in our paper is analogous to the azimuthal state on an annulus, with tangent boundary conditions on both circular boundaries \cite{duffy1, duffy2, duffy3, palffy}. Whilst our results in the OF case, both with strong and weak anchoring, are partially captured by the results in \cite{pan, duffy1, duffy2, duffy3, palffy}, our method of proof is different. We compute the second variation of the anisotropic OF energy and study the resulting eigenvalue problem directly. We compute stability curves in terms of the elastic anisotropy $\delta$, the annular aspect ratio $b$, the anchoring strength, $\alpha$, on both boundaries and the order of the azimuthal perturbations, $k$. Such stability diagrams can provide useful insight into stabilization and de-stabilization effects in simple geometries, i.e.~they can quantify the response of the defect-free state to different types of azimuthal perturbations. For $k=0$, our results reduce to the previously reported results in \cite{pan, duffy1, duffy2, duffy3, palffy}. Finally, in \cite{helein}, the authors demonstrate that the planar radial state loses stability, when the ratio of the inner radius to the outer radius is smaller than a critical value, and the planar radial state escapes into the `third' dimension. In \cite{virga2}, the authors  study nematic samples confined between two co-axial cylinders, with emphasis on higher-dimensional biaxial effects, which are outside the scope of the present paper.

In \cite{mottramdavidson, lewis1}, the authors study nematic equilibria confined to shallow square and rectangular wells, subject to strong tangent boundary conditions on the edges. They construct analytic approximations for the experimentally observed `diagonal' and `rotated' states and obtain explicit expressions for their respective energies. The diagonal state always has lower energy than the rotated states. We take this work further by constructing explicit solutions of the Laplace's equation on an annular sector of angle $2\pi/N$ (with $N \in \mathbb{N}$) with tangent boundary conditions. These `sector' solutions mimic generalized `rotated' and `diagonal' nematic equilibria on 2D annuli, featuring a total of $2N$ defects pinned to the inner and outer boundaries. We compute the respective energies in the OF framework and find that there are interesting energetic cross-overs, in terms of the number of defects. Within a sector, we find that the rotated states typically have lower energies than the corresponding diagonal states for small values of $N$. As $N$ increases, the sector approaches a rectangle and the generalized diagonal state becomes energetically preferable, consistent with the energetic trends in a rectangle. 
States with defects on the boundary of an annular well are considered in \cite{Viral_Nematics}. The authors work within the one constant Oseen-Frank framework with weak anchoring and use an expression for the director based on a linear combination of defect-type solutions located at equally spaced locations on the boundary. This expression, although not a solution of the Euler-Lagrange equations, allows the authors to identify possible director profiles and obtain estimates to the surface and elastic energies. These energies demonstrate that states with boundary defects may have lower energy than the defect-free state in certain parameter regimes by including weak anchoring. In this paper, we demonstrate that in the strong anchoring regime, this is only possible for liquid crystals with large elastic anisotropy.



Equilibria with boundary defects can be interesting in physical situations, particularly with Neumann boundary conditions or when we simply specify topological degrees on the boundaries, as opposed to Dirichlet conditions. For example, in \cite{beryland}, the authors study 2D vector fields, $\uvec = (u_1, u_2)$, on a multiply-connected 2D domain with prescribed topological degrees on the outer boundary and on the inner boundaries enclosing the `holes'. The authors study minimizers of the Ginzburg-Landau functional on such domains and find that the infimum energy is not attained. Minimizing sequences develop vortices or boundary defects in certain asymptotic limits. Our analysis of OF equilibria with  boundary defects could be instructive for such model LC problems, with topological degree boundary conditions; in fact, they give quantitative information about minimizing sequences. Equally importantly, our analysis is directly relevant to recent experimental work on colloidal samples in shallow annular wells wherein the authors observe states with  defects pinned to the lateral surfaces \cite{Oli_thesis, Jose_thesis}, such as those shown in Figure \ref{fig:Intro}. The relative energies of the generalized `diagonal' and `rotated' states give qualitative insight into the relative observational frequencies of the experimental states that exhibit boundary defects.

We complement our analysis in the OF case with work within the more general and powerful Landau-de Gennes (LdG) theory for nematic liquid crystals.
In \cite{mironescu}, the author studies the generic $+1$-degree vortex within the Ginzburg-Landau (GL) theory for superconductivity and derives a powerful stability result for the $+1$-degree vortex on a 2D disk (that contains the origin), with Dirichlet radial boundary conditions. In 2D, the LdG theory for nematic liquid crystals reduces to the GL theory and the nematic state is fully described by a scalar order parameter and a two-dimensional vector field, referred to as `director' in the literature \cite{maj1}. We define the nematic defect-free state in the LdG framework, in terms of a minimizer of an appropriately defined energy functional, and adapt Mironescu's proof to a 2D annulus. We address technical differences due to the change in geometry and demonstrate local stability in the low temperature limit, and derive a stability criteria in terms of temperature and geometry.

The paper is organized as follows. In Section~\ref{sec:theory}, we review the OF and LdG theory for nematic LCs. In Section~\ref{sec:strong} and \ref{sec:weak}, we focus on the OF theory and compute stability criteria for the defect-free state  as a function of $\delta$, $b$ and $\alpha$ as introduced above. In Section~\ref{sec:defects}, we construct generalized `diagonal' and `rotated' states in an annulus with an arbitrary number of boundary defects and compute the corresponding energies. Finally, in Section~\ref{sec:mironescu}, we compute criteria for the local minimality of the LdG defect-free state on a 2D annulus as described above and conclude in Section~\ref{sec:conclusion} with future perspectives.

\section{Theory and Modelling}
\label{sec:theory}

Nematic liquid crystals are complex liquids wherein the constituent rod-like molecules move freely as in a conventional liquid and tend to align along certain  preferred directions \cite{dg,virga2, newtonmottram}. The Landau-de Gennes (LdG) theory is one of the most general continuum theories for nematic LCs to date \cite{dg, newtonmottram}. The LdG theory describes the nematic state by a macroscopic order parameter, known as the LdG $\Qvec$-tensor, that is a macroscopic measure of the LC anisotropy. For three-dimensional problems, the LdG $\Qvec$-tensor is a symmetric, traceless $3 \times 3$ matrix with five degrees of freedom \cite{dg,newtonmottram}. A nematic phase is said to be (i) isotropic if $\Qvec=0$, (ii) uniaxial if $\Qvec$ has two degenerate non-zero eigenvalues and (iii) biaxial if $\Qvec$ has three distinct eigenvalues. The LdG theory is a variational theory and a prototypical LdG energy functional is of the form
\begin{equation}
\label{eq:1}
I[\Qvec]:= \int_{\Omega} w(\Qvec, \grad \Qvec) + f_B (\Qvec) ~ d\Omega  \quad \textrm{where $f_B(\Qvec) = \frac{A}{2}\textrm{tr}\Qvec^2 - \frac{B}{3}\textrm{tr}\Qvec^3 + \frac{C}{4}\left(\textrm{tr} \Qvec^2 \right)^2$,}
\end{equation}
$w(\Qvec, \grad \Qvec)$ is an elastic energy density that penalizes spatial inhomogeneities and $f_B$ is a bulk potential that drives the nematic-isotropic phase transition as a function of the temperature or concentration. In Equation (\ref{eq:1}), $A<0$ and is proportional to the rescaled temperature and $B, C>0$ are positive material-dependent constants. Whilst studying 2D nematic equilibria on 2D domains, the LdG $\Qvec$-tensor reduces to a symmetric, traceless $2\times 2$ matrix of the form 
\begin{equation}
\label{eq:2}
\Qvec = s(x,y) \left( \nvec \otimes \nvec - \frac{I_2}{2} \right),
\end{equation}
where $\nvec$ is a two-dimensional unit-vector and $I_2$ is the identity matrix in 2D \cite{maj1}. It is clear from (\ref{eq:2}) that a 2D $\Qvec$-tensor only has two degrees of freedom: the real-valued scalar order parameter, $s(x,y)$, that is a measure of the degree of the orientational order and the unit-vector field, $\nvec \in S^1$, that represents the distinguished direction of alignment of the nematic molecules. For $\Qvec$ as in (\ref{eq:2}), $\textrm{tr}\Qvec^3 = 0$ and hence, the LdG energy reduces to
 \begin{equation}
\label{eq:3}
I[\Qvec]:= \int_{\Omega} w(\Qvec, \grad \Qvec) + \frac{A}{2}\textrm{tr}\Qvec^2  + \frac{C}{4}\left(\textrm{tr} \Qvec^2 \right)^2 d \Omega .
\end{equation} We work with low temperatures, described by $A<0$, and the one-constant elastic energy density, $w(\Qvec, \grad\Qvec) = \frac{L}{2} |\grad \Qvec|^2$. We introduce the scaling $\bar{\Qvec} = \frac{\Qvec_{ij}}{\sqrt{|A|/ C}}$, so that the LdG energy in (\ref{eq:3}) reduces to the Ginzburg-Landau functional in superconductivity \cite{maj1,bbh},
\begin{equation}\label{eq:3b}
I[\Qvec] = \int_{\Omega} |\grad \Qvec|^2 + \frac{\vert A \vert}{4L} \left( \vert \Qvec \vert^2 - 1\right)^2 d \Omega .
\end{equation}
 As in any problem in the calculus of variations, the problem of studying nematic equilibria is mathematically equivalent to a  study of local and global energy minimizers.
In particular, we obtain an explicit relation between geometry and temperature that guarantees local stability of the defect-free state.

The Oseen-Frank (OF) theory is a simpler continuum theory for nematic liquid crystals restricted to uniaxial nematic phases with a constant scalar order parameter \cite{dg,virga}. In this case, the macroscopic order parameter is simply a unit-vector field, $\nvec$, that defines the unique direction of molecular alignment. For a 2D problem, the OF energy functional is given by
\begin{equation}
\label{eq:4}
E[\nvec]: = \iint_{\Omega} \frac{K_1}{2}\left( \nabla \cdot \nvec \right)^2 + \frac{K_3}{2}\left(\nvec \times \left(\grad \times \nvec \right) \right)^2~d\Omega ,
\end{equation}
where $\Omega \subset \Rr^2$ is a 2D domain and the elastic constants $K_1$ and $K_3$ are associated with splay and bend director deformations \cite{dg}.
For 2D domains, $\nvec$ can be conveniently written as
$
\nvec = \left( \cos \theta, \sin \theta, 0 \right),
$
where $\theta: \Rr^2 \to \Rr$ is a function of the planar polar coordinates, $(r, \phi)$, and the OF energy (\ref{eq:4}) reduces to
\begin{eqnarray}
\label{eq:6}
E[\theta] :=\iint_\Omega   \frac{K_1}{2} \left( \cos(\theta - \phi) \frac{\theta_\phi}{r} - \sin(\theta - \phi) \theta_r \right)^2  + \frac{K_3}{2} \left( \sin(\theta - \phi) \frac{\theta_\phi}{r} + \cos(\theta - \phi) \theta_r \right)^2 ~d \Omega .
\end{eqnarray}
We define $\delta := 1 - \frac{K_1}{K_3}$ to be the measure of elastic anisotropy. The corresponding Euler-Lagrange equations are:
\begin{eqnarray}\label{eq:7}
&& \nabla^2 \theta \left( 1 - \frac{\delta}{2}\right) +\delta \left( \frac{1}{2} \sin (2\theta - 2\phi) \left( \frac{2 \theta_{r\phi}}{r} + \frac{\theta_\phi^2}{r^2} - \theta_r^2 - \frac{2\theta_\phi}{r^2} \right)\right.  \nonumber \\
&&  \qquad \qquad \qquad \qquad + \left. \frac{1}{2} \cos(2\theta - 2 \phi) \left( \theta_{rr} - \frac{\theta_{\phi \phi}}{r^2} - \frac{\theta_r}{r} + \frac{2\theta_r \theta_\phi}{r} \right) \right) =0 ,
\end{eqnarray} subject to appropriate boundary conditions for $\theta$ on $\partial \Omega$.

In Sections~\ref{sec:strong}, \ref{sec:weak} and \ref{sec:mironescu}, we focus on the defect-free state on a 2D annulus with tangent boundary conditions. The rescaled 2D annulus is defined by
\begin{equation}
\label{eq:8}
\Omega = \left\{ (r, \phi) \in \Rr^2; \quad  b \leq r \leq 1, \quad 0\leq \phi < 2 \pi \right\},
\end{equation} where the rescaled radius $b := \frac{R_{\text{inner}}}{R_{\text{outer}}}$ is the ratio of the inner and outer radii. In Section~\ref{sec:defects}, we consider the reduced domain
\begin{equation}
\Omega_N =\left\{ (r, \phi) \in \Rr^2; \quad b \leq r \leq 1, \quad 0\leq \phi < \frac{2 \pi}{N} \right\},
\end{equation}
where $N  \in  \mathbb{N}$, which allows us to consider more complex equilibria with defects.
 The tangent boundary conditions require that
\begin{equation}
\label{eq:9}
\theta(1, \phi)  = \theta(b, \phi) = \phi \pm \frac{\pi}{2}.
\end{equation}
Equation~(\ref{eq:9}) is the Dirichlet version or the strong version of the tangent boundary conditions. The tangent boundary conditions can be weakly implemented too, as will be illustrated in Section~\ref{sec:weak}.
The OF defect-free state is simply defined by
\begin{equation}
\label{eq:10}
\theta^*(r, \phi) = \phi + \frac{\pi}{2}, \qquad (r, \phi) \in \Omega,
\end{equation} with no radial variations. A simple computation shows that $E[\theta^*] = \pi K_3 \log \left(\frac{1}{b} \right)
,$ so that the OF energy diverges in the limit $b \to 0$. In the LdG framework, by analogy with similar work on degree $+1$-vortices in the GL theory \cite{bbh}, we define the LdG defect-free state to be
\begin{equation}
\label{eq:12}
\Qvec^*(r, \phi) = s\left(r \right) \left(\nvec^* \otimes \nvec^* - \frac{I_2}{2} \right)
\end{equation}
where
$$ \nvec^* = \left( \cos \theta^*, \sin \theta^*, 0 \right) $$
and $s$ is the unknown scalar order parameter, independent of $\phi$.
In what follows, we study the defect-free state in both theoretical frameworks and make comparisons to competing states with boundary defects.

\section{The OF Defect-Free State with Strong Anchoring}
\label{sec:strong}

We firstly consider the OF defect-free state in (\ref{eq:10}) with Dirichlet boundary conditions. Our result is similar to the stability estimate in \cite{pan, palffy} but our method of proof is different. We consider the second variation of the OF energy about $\theta^*$. The strict positivity of the second variation is a sufficient criterion for local stability, whereas $\theta^*$ is unstable if the second variation is negative for some admissible perturbation \cite{hestenes}.

We perturb $\theta^*$ by $\epsilon \eta$, where
$ \eta(b, \phi) = \eta(1, \phi) = 0$
in accordance with the Dirichlet conditions in (\ref{eq:9}). A straightforward computation shows that the second variation is given by
\begin{eqnarray}
\label{eq:15}
\delta^2 E: = \frac{\partial^2 E[\theta^* + \epsilon \eta]}{\partial \epsilon^2} \Bigg \vert_{\epsilon=0} &=& \iint_\Omega (K_1 - K_3) \left(\frac{\eta}{r} + \eta_r \right)^2 + K_3 (\nabla \eta)^2 ~d \Omega \\
&=& K_3 \iint_\Omega (\nabla \eta)^2 - \delta \left(\frac{\eta}{r} + \eta_r \right)^2 ~d \Omega
\end{eqnarray}
and it follows immediately that $\delta^2 E > 0$ for $\delta \leq 0$ and any non-trivial $\eta$.
\paragraph{Proposition 1:}
For any admissible $\eta(r)$, there exists $\delta_\eta < 1$ such that $\delta^2 E < 0$ for $\delta_\eta < \delta \leq 1. $

\begin{proof}
We define the functional $E_2$,
\begin{equation}
E_2 (\eta,\delta,b) := \iint_\Omega (\nabla \eta)^2 - \delta \left(\frac{\eta}{r} + \eta_r \right)^2 d A.
\end{equation}
For $\eta$ independent of $\phi$, $E_2$ simplifies to
\begin{equation}
E_2 (\eta(r),\delta,b) = 2\pi \int^1_b \left\lbrace  (\eta_r)^2  - \delta \left(\frac{\eta}{r} + \eta_r \right)^2 \right\rbrace r d r.
\end{equation}
For a fixed $b$ and $\eta(r)$, $E_2$ is a continuous function of $\delta$ which is bounded below by $E_2 (\eta(r),1,b)$.
Using the boundary conditions $\eta(1) = \eta(b) = 0$, this lower bound simplifies to
\begin{eqnarray}
\label{eq:16}
E_2 (\eta(r),1,b) &=&-2\pi \int^1_b \frac{\eta^2}{r} + 2\eta \eta_r dr = -2\pi \int^1_b \frac{\eta^2}{r} d r
\end{eqnarray}
and hence
 $
E_2 (\eta(r),1,b) < 0
$
for all $\eta$ and any $0 < b < 1$. We use the continuity of $E_2$ with respect to $\delta$ to conclude that for all $\eta(r)$, there exists $\delta_\eta < 1$ such that $\delta^2 E < 0$ for all $\delta_\eta < \delta \leq 1$.
\end{proof}
The preceding observation demonstrates that the defect-free state loses stability on a 2D annulus, for certain choices of $(\delta, b)$. We make this observation more rigorous by studying the Sturm-Liouville problem associated with the minimization of $\delta^2 E$ as shown below. The corresponding Euler-Lagrange equation is:
\begin{equation} \label{eq:17}
\frac{\delta \eta}{r^2} + \frac{\eta_{\phi \phi}}{r^2} + (1-\delta) \left(\eta_{rr} + \frac{\eta_r}{r} \right) =0,
\end{equation}
where $\eta = 0$ on $r = b,1$ and $\eta$ is $2\pi$-periodic.
Without loss of generality, we can consider separable solutions of the form
$
\eta = \sum_{k} e^{ik\phi}f_k(r),
$
where $k\in \mathbb{Z}$.
The function $f_k(r)$ is a solution of
\begin{equation} \label{eq:19}
r f_k'(r)+ r^2 f_k''(r) + \frac{\delta - k^2}{1-\delta} f_k(r) = 0,
\end{equation}
with boundary conditions $f_k(b) = f_k(1) = 0$.
Equation (\ref{eq:19}) is an Euler ordinary differential equation and we take $f_k$ to be of the form $f_k = r^m$, where $m$ and $k$ are related by
\begin{equation} \label{eq:20}
m^2 = \frac{ k^2 - \delta}{1-\delta}.
\end{equation}
For $k=0$, the general solution of (\ref{eq:19}) is given by
\begin{equation} \label{eq:21}
f_0(r) = A \sin \left( \sqrt{\frac{\delta}{1- \delta}} \log (r) \right) + B \cos \left( \sqrt{\frac{\delta}{1- \delta}} \log (r) \right)
\end{equation}
and applying the boundary conditions $f_0(1)= f_0(b) =0$, we find that
\begin{equation} \label{eq:22}
f_{0,n}(r) = A\sin \left( \pi n \frac{\log (r)}{\log (b)} \right) \quad \textrm{where $\delta_n = \frac{\pi^2 n^2}{\pi^2 n^2 + \log(b)^2}$. }
\end{equation}
For $k \geq 1$, the boundary conditions on $f_k$ naturally rule out any non-trivial solutions. The case $k=0$ is a typical Sturm-Liouville problem \cite{sturm-liouville} with eigenfunctions $\left\{f_{0,n} \right\}$ and corresponding eigenvalues, $\left\{ \delta_n \right\}$. 
 Similarly, substituting
 $\eta  = A \sin \left( \pi  \frac{\log (r)}{\log (b)} \right)$ into (\ref{eq:15}), we find that the corresponding second variation is negative, $\delta^2 E < 0$ for $\delta > \delta_1$. Therefore, the defect-free state loses stability along the curve, $\delta = \delta_1\left( b \right)$ defined by
\begin{equation}
\label{eq:24}
\delta_1 = \frac{\pi^2 }{\pi^2  + \log(b)^2}.
\end{equation}
This stability criterion (\ref{eq:24}) is known in the literature \cite{pan,duffy3}.

We next establish that the defect-free state undergoes a supercritical pitchfork bifurcation at $\delta = \delta_1$ above. We let $\delta = \delta_1 + \epsilon^2 \delta_2$ and consider perturbations of the form
\begin{equation} \label{eq:theta}
\theta =\theta^* + \epsilon \eta_1 + \epsilon^2 \eta_2 + \epsilon^3 \eta_3 + \dots
\end{equation} where
the functions $\eta_i$ vanish on $r = b,1$ and are $2\pi$-periodic in $\phi$. We substitute (\ref{eq:theta}) into the Euler-Lagrange equation (\ref{eq:7}) and compare terms of order $\eps$. We define the operator $\mathcal{L}$ to be
\begin{equation}
\label{eq:L}
\mathcal{L}\eta:= \frac{\delta_1 \eta}{r^2} + \frac{\eta_{\phi \phi}}{r^2} + (1 - \delta_1) \left( \eta_{rr} + \frac{\eta_{r}}{r} \right).
\end{equation}
A straightforward computation shows that $\mathcal {L}\eta_1  = \mathcal{L}\eta_2 = 0 $. Recalling (\ref{eq:17}), we deduce that $\eta_1$ and $\eta_2$ are proportional to $ \sin\left(\pi \frac{\log(r)}{\log(b)} \right)$. Comparing terms of order $\eps^3$, we find
\begin{equation} \label{eq:L3}
\mathcal{L} \eta_3 = R(r) ,
\end{equation}
where the function $R$ is defined to be
\begin{align}\label{eq:L4}
R(r) :=&  \frac{(\ln^2(b) + \pi^2)A(\delta_1 A^2 - 2\delta_2)}{2r^2 \ln^2 (b)} \sin\left(\pi\frac{\ln(r)}{\ln(b)} \right) \nonumber \\ &-\frac{\delta_1 A^3 (3\pi^2 + \ln^2 (b))}{6r^2 \ln^2 (b)}  \sin\left(3\pi\frac{\ln(r)}{\ln(b)}\right).
\end{align}
 From the Fredholm Alternative Theorem, the differential equation (\ref{eq:L3}) admits a solution if and only if the following compatibility condition holds \cite{sturm-liouville}:
\begin{equation}
\int^1_b  \sin\left(\pi \frac{\log(r)}{\log(b)} \right) R(r)r d r =  0,
\end{equation}
This condition reduces to
\begin{equation} \label{eq:L5}
A ( A^2 \delta_1 - 2\delta_2) =0,
\end{equation}
and therefore the amplitude, $A$, is only defined for $\delta_2>0$ and we have a supercritical pitchfork bifurcation at $\delta=\delta_1$.

Finally, we explicitly construct a solution of the Euler-Lagrange equations (\ref{eq:7}) that gives us insight into how the defect-free solution deforms for $\delta > \delta_1$ and the structure of the corresponding global energy minimizers for $\delta > \delta_1$. We seek an exact solution of the form
\begin{equation}
\label{eq:25}
\theta = \theta^* + U(\delta, r),
\end{equation}
where $U$ satisfies $U(\delta, 1) = U(\delta, b) =0$. We consider functions with at most one critical point, solutions with more critical points exist but have a higher energy. We define $t = -\log(r)$, substitute (\ref{eq:25}) into (\ref{eq:7}) and then $U$ satisfies
\begin{equation}\label{eq:27}
(\delta \cos^2 (U(t)) - 1) \frac{\text{d}^2 U}{\text{d} t^2} - \frac{\delta}{2} \sin(2U(t)) \left( \left( \frac{\text{d} U}{\text{d} t}\right)^2 + 1 \right) =0,
\end{equation}
subject to $U(\delta,0)  = 0$ and $U(\delta, \log(1/b))  = 0$.

We can define $U(\delta,t)$ implicitly by
\begin{equation}
\int^{U(\delta,t)}_0 \sqrt{\frac{1 - \delta \cos^2 u }{ \delta \cos^2 u - \delta \cos^2 U_{0}}} d u = t,
\end{equation}
where $U_0$ is the maximum value of $U(\delta, t)$.

\begin{figure}[ht]
\centering
\includegraphics[scale=0.8]{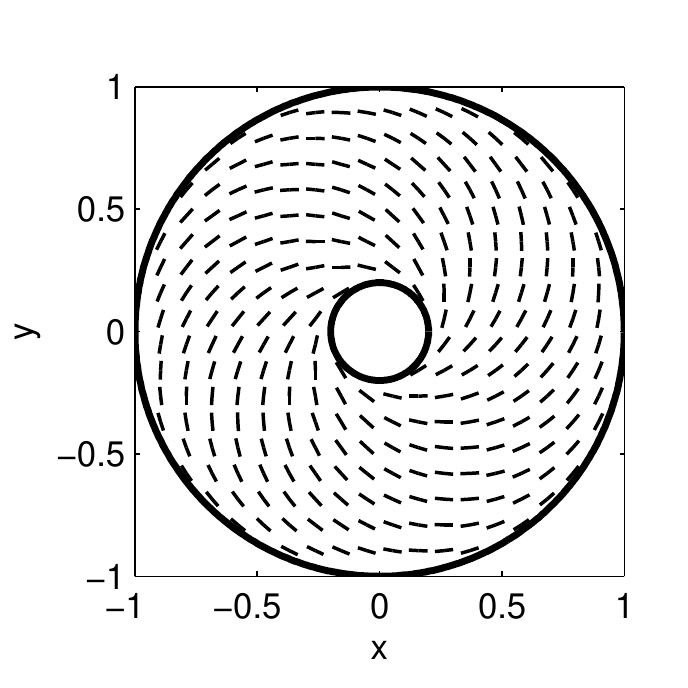}
\label{fig:strong_anchoring}
\caption{The director field for $\theta = \theta^* + U$, with $b = 0.2$ and $\delta =0.95$.}
\end{figure}

 In Figure~\ref{fig:strong_anchoring}, we plot the spiral-like solution (\ref{eq:25}) for $b=0.2$ and $\delta=0.95$. Numerical computations of the energy show that this spiral-like solution has lower energy than the defect-free solution in the parameter regime $\delta > \delta_1$.

When $\delta = 1$, we can solve (\ref{eq:27}) exactly to find
\begin{equation}\label{eq:29}
U(1,t) = \arccos \left( \frac{b}{b+1} \exp(t) + \frac{1}{b+1} \exp(-t)\right)
\end{equation}
with associated energy, $E = 2\pi\left( \frac{1 -b}{1+ b}\right)$.
Given that (\ref{eq:29}) is an exact solution, we can study its stability in terms of the second variation;
\begin{equation}
\delta^2 E: = \iint_\Omega g^2 \eta_\phi^2 + (1-g^2) \eta_t^2 + \eta^2 \left( \frac{2 g^4 - 3g^2 + g'^2 + 1 + 2 g g'' - 2 g^3 g''}{1 - g^2}\right) d \Omega ,
\end{equation}
where $g(t) =\frac{b}{b+1} \exp(t) + \frac{1}{b+1} \exp(-t)$.
The function $g \leq 1$ for $t \in (0,-\log(b))$ and the coefficient of $\eta^2$ simplifies to
\begin{equation}
\frac{2b - b^2 -1}{b^2 e^{2t} + e^{-2t} - b^2 - 1} \geq 1.
\end{equation}
Therefore, the second variation is positive for all admissible perturbations and we have established local stability of the spiral-like solution (\ref{eq:25}) at $\delta=1$. In summary, we have constructed an exact solution of the Euler-Lagrange equations, (\ref{eq:7}), which is a radial variant of the defect-free state and is locally stable for $\delta$ close to unity.

\section{The OF Defect-Free State with Weak Anchoring}
\label{sec:weak}

We add a surface energy to the Oseen-Frank energy in (\ref{eq:4}), that enforces preferential planar anchoring. Our approach is similar to the work in \cite{palffy} and our main contribution is the stability diagram in Figure~\ref{fig:weak} that is particularly useful in quantifying the interplay between anisotropy, geometry and surface effects in the stability of the defect-free state. We employ the Rapini-Papoular surface energy \cite{RapiniPapoular}
\begin{equation} \label{eq:s1}
E_S = \frac{1}{2} \int_{\partial \Omega}  W R_{\text{outer}} \sin^2 \left(\theta - \phi - \frac{\pi}{2}\right) d s,
\end{equation}
where $\partial \Omega$ consists of two concentric circles with radii $r=b$ and $r=1$, $W$ is an anchoring coefficient 
and $ds$ is an arc-length element along $\partial \Omega$. We define $\alpha := \frac{W  R_{\text{outer}}}{K_3} = \frac{R_{\text{outer}}}{\xi}$ to be the ratio of the outer radius to the extrapolation length $\xi$, \cite{lewis1}, and the three key parameters are $\delta$, $b$ and $\alpha$. The corresponding Euler-Lagrange equation is (\ref{eq:7})
with boundary conditions
\begin{eqnarray}\label{eq:s3}
\left(\frac{2 - \delta}{2}\right) \theta_r + \frac{\delta}{2} \left( \frac{\theta_\phi}{r} \sin(2\theta - 2\phi) + \theta_r \cos(2\theta - 2\phi)\right) +\frac{\alpha}{2} \sin(2\theta - 2\theta_0) =0,
\end{eqnarray}
on $r=1$ and
\begin{eqnarray}\label{eq:s4}
-\left(\frac{2 - \delta}{2}\right) \theta_r - \frac{\delta}{2} \left( \frac{\theta_\phi}{r} \sin(2\theta - 2\phi) + \theta_r \cos(2\theta - 2\phi)\right) +  \frac{\alpha}{2} \sin(2\theta- 2\theta_0) =0,
\end{eqnarray}
on $r=b$. It is straightforward to check that the OF defect-free state, $\theta^*$, defined in (\ref{eq:10}), is a solution of (\ref{eq:7}) with the above boundary conditions. By analogy with our work in Section~\ref{sec:strong}, we compute the second variation as shown below:
\begin{equation}\label{eq:secondweak}
 \delta^2 E[\theta^*]: = \iint_\Omega  (\nabla \eta)^2 - \delta \left(\frac{\eta}{r} + \eta_r \right)^2 d \Omega + \alpha  \int_{\partial \Omega} \eta^2 d s 
\end{equation}
where $\eta$ is a perturbation about $\theta^*$. 
Without loss of generality, we assume $
\eta = \sum_{k} f_k(r) e^{ik\phi},
$ where we refer to $k$ as being the azimuthal order of the perturbation. For each $k$, the optimal $f_k$ is a solution of
\begin{equation} \label{eq:s6}
rf_{k}'(r) + r^2 f_{k}''(r) + \frac{\delta - k^2}{ 1-\delta}f_k (r)  =0,
\end{equation}
with boundary conditions
\begin{eqnarray}\label{eq:s7a}
f_{k}'(r) =\frac{\delta - \alpha}{1-\delta} f_k(r) \quad \text{on $r = 1$,}\\
f_{k}'(r)= \frac{\alpha + \frac{\delta}{b}}{1-\delta} f_k(r) \quad \text{on $r = b$.} \label{eq:s7b}
\end{eqnarray}
Standard computations, following a change of variable $t = \log \left(\frac{1}{r}\right)$, show that $f_k(t)$ is given by
\begin{align}
f_0(t) &= A \sin \left( \sqrt{\frac{ \delta}{1-\delta}} t\right) + B \cos \left( \sqrt{\frac{\delta}{1-\delta}}t\right)  \qquad \text{for $k=0$} \\
f_k(t) &= A \sinh \left( \sqrt{\frac{k^2 - \delta}{1-\delta}} t\right) + B \cosh \left( \sqrt{\frac{k^2 - \delta}{1-\delta}}t\right) \qquad \text{for $k\geq 1$}
\end{align}
and the boundary conditions require the following `compatibility condition' between $\alpha, \delta, b$ and $k$:

\begin{align} 
&\tan  \left( \sqrt{\frac{ \delta}{1-\delta}} \ln  \left( \frac{1}{b} \right)\right) + \frac{\alpha ( 1 + b) \sqrt{\delta ( 1 - \delta)}}{\alpha \delta - \alpha b \delta + \alpha^2 b - \delta} = 0 \qquad \text{for $k =0$} \label{eq:s11}  \\
&\tanh \left(\sqrt {{
\frac {{k}^{2}-\delta}{1-\delta}}}\ln  \left( \frac{1}{b} \right) \right) + \frac{(1-\delta) \alpha (1+b)}{\delta \alpha + \alpha^2 b - \alpha b \delta - \delta + k^2 ( 1-\delta)} \sqrt {
\frac {{k}^{2}-\delta}{1-\delta}} =0 \qquad \text{for $k \geq 1$} . \label{eqn:s11.1}
\end{align}
For a given $k \in \mathbb{N}$, there will typically be multiple solutions $\left\{\delta_{1,k}, \delta_{2,k}, \ldots \right\}$ of the compatibility relation above. We are interested in the smallest solutions, $\delta_{1,k}$, whilst defining stability curves in the $\left(\delta, \alpha \right)$-plane for different values of $k$ and $b$.

\subsection{The case $k=0$}
\begin{figure}[ht]
\centering
\includegraphics[scale=0.8]{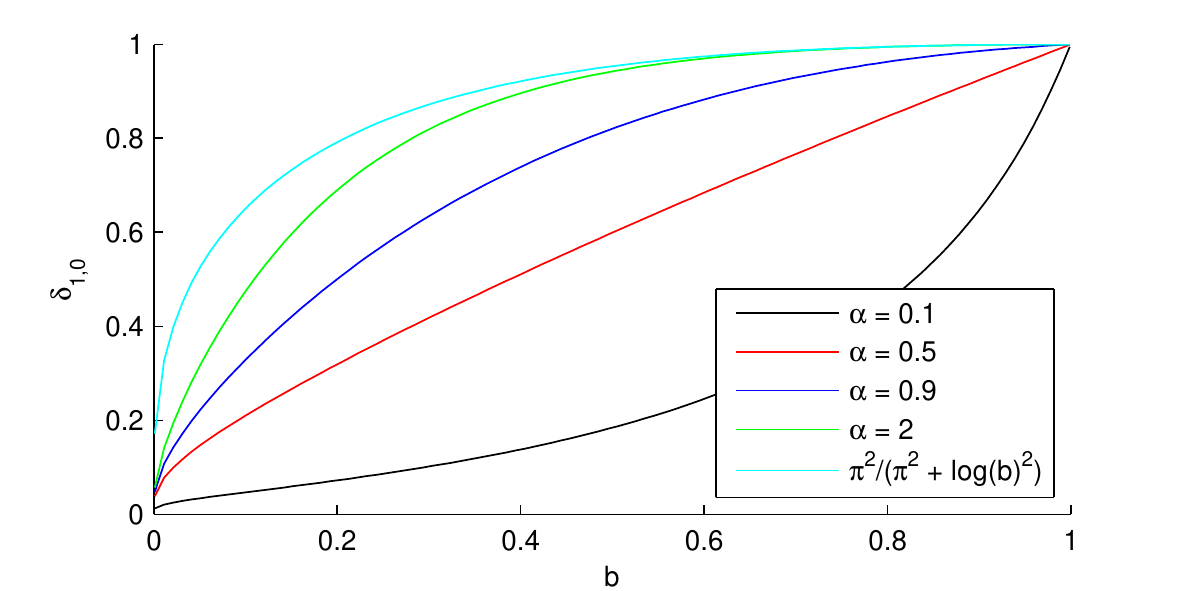}
\caption{The first eigenvalue $\delta_{1,0}$ for fixed $\alpha$ and the strong anchoring limit.}
\label{fig:weak1}
\end{figure}

The case $k=0$ is contained in the work of \cite{palffy}. However, our method of proof is different and of independent interest. We first point out that for $k=0$, the boundary-value  problem (\ref{eq:s6}) - (\ref{eq:s7b}) is a Sturm-Liouville problem with a complete set of eigenfunctions and corresponding eigenvalues, $\left\{ \delta_{1,0}, \delta_{2,0}, \ldots, \delta_{n,0} \right\}$ \cite{sturm-liouville}. The  eigenvalues, $\left\{\delta_{n,0}\right\}$, are solutions of the compatibility condition \begin{equation}
\tan  \left( \sqrt{\frac{ \delta}{1-\delta}} \ln  \left( \frac{1}{b} \right)\right) + \frac{\alpha ( 1 + b) \sqrt{\delta ( 1 - \delta)}}{\alpha \delta - \alpha b \delta + \alpha^2 b - \delta} = 0,\label{eq:s12}
\end{equation} which has an infinite number of solutions for any fixed triplet $\left(\alpha, \delta, b \right)$. The corresponding eigenfunctions are given by
\begin{equation}\label{eq:s13}
\eta_{n,0} (r,\phi) = A \sin \left( \sqrt{\frac{ \delta_{n, 0}}{1-\delta_{n, 0}}} \ln  \left( \frac{1}{r} \right)\right) + A \frac{\sqrt{\delta_{n, 0}(1-\delta_{n, 0})}}{\alpha - \delta_{n, 0}} \cos \left( \sqrt{\frac{ \delta_{n, 0}}{1-\delta_{n, 0}}}\ln  \left( \frac{1}{r} \right)\right).
\end{equation} In Figure~\ref{fig:weak1}, we compute $\delta_{1,0}$ as a function of $b$, for different values of $\alpha$. As $\alpha \rightarrow \infty$, the compatibility condition becomes $\tan  \left( \sqrt{\frac{ \delta}{1-\delta}} \ln  \left( \frac{1}{b} \right)\right) \rightarrow 0$ and $\delta_{1,0} \rightarrow \frac{\pi^2}{\pi^2 + \ln(b)^2}$, recovering the strong anchoring result.

We further investigate the loss of stability of $\theta^*$ at $\delta=\delta_{1,0}$. Let $\theta = \theta^* + \epsilon \eta_{1} + \epsilon^2 \eta_2 + \epsilon^3 \eta_3$ and $\delta = \delta_{1,0} + \epsilon^2 \delta_2$ as before. Then $\eta_1$ and $\eta_2$ are p roportional to Equation (\ref{eq:s13}) with $n=1$. The corresponding compatibility condition for $\eta_3$ to exist is
\begin{align}
\int^1_{b} \eta_{3h}  R(r) r d r &= r\eta_{3h}\left(\delta_2 \left(\frac{\eta_{1}}{r} + \eta_{{1}}'\right) + \frac{2}{3} \eta_1^3 \left( \alpha - \frac{\delta_{1,0}}{r}\right) - \delta_{1,0} \eta_1^2 \eta_{{1}}'\right)\bigg
\vert_{r=1} \nonumber \\
& - r\eta_{3h} \left(\delta_2 \left(\frac{\eta_{1}}{r} + \eta_{{1}}'\right) - \frac{2}{3} \eta_1^3 \left(\alpha + \frac{\delta_{1,0}}{b} \right)- \delta_{1,0} \eta_1^2 \eta_{{1}}' \right) \bigg \vert_{r=b},
\end{align}
where $\eta_{3h}$ is the solution of the homogeneous problem as given below:
\begin{equation}
\eta_{3h} (r,\phi) =  \sin \left( \sqrt{\frac{ \delta_{1,0}}{1-\delta_{1,0}}} \ln  \left( \frac{1}{r} \right)\right) +  \frac{\sqrt{\delta_{1,0}(1-\delta_{1,0})}}{\alpha - \delta_{1,0}} \cos \left( \sqrt{\frac{ \delta_{1,0}}{1-\delta_{1,0}}}\ln  \left( \frac{1}{r} \right)\right),
\end{equation}
and
\begin{equation}
R(r) \equiv \delta_2 \left( \frac{\eta_{{1}}'}{r} + \eta_{{1}}'' - \frac{\eta_{1}}{r^2} \right) - \delta_{1,0} \left( \eta_{1} \eta_{{1}}'^2  - \frac{2\eta_{{1}}^3}{3r^2} + \eta^2_{1} \eta_{{1}}''  + \frac{\eta_{1}^2 \eta_{{1}}'}{r}\right).
\end{equation}
The compatibility condition yields
\begin{equation}
A^3 E_3 - \delta_2 A E_1 = 0,
\end{equation}
where $E_3$ and $E_1$ are functions of $\delta_{1,0}$, $b$ and $\alpha$. We numerically can evaluate $E_1$ and $E_3$ at $\delta =  \delta_{1,0} (\alpha,b)$ to find that $E_1$ and $E_3$ are both positive.
Hence, we have a supercritical pitchfork bifurcation at $\delta=\delta_{1,0}$ in the weak anchoring framework.

\subsection{The cases $k\geq 1$}

For $k=1$, the compatibility condition (\ref{eqn:s11.1}) simplifies to \begin{equation}
-\frac{(1+ b) (b^2 \alpha^2 - \alpha b^2 + 2\delta\alpha b - 2\delta b + b - \alpha^2 b - \alpha + 2\delta -1)}{(b^2 + 1 )(\delta \alpha + \alpha^2 b - \delta \alpha b - 2\delta + 1}=0,
\end{equation}
with a unique solution
\begin{equation}
\delta_{1,1}\left(\alpha, b \right) = \frac{1}{2}\frac{\alpha b^2 + \alpha^2 b + \alpha + 1 - b^2 \alpha^2 - b}{\alpha b - b+1}.
\end{equation}
We note that $\delta_{1,1} \in (0,1)$ iff $0 <\alpha <1$ and that $\delta_{1,1} > 0.5$ for all $b$.

\begin{figure}[!ht]
        \begin{center}
        \subfigure[$b=0.1$]    {
                        \label{Stability_curveb=01}
                        \includegraphics[scale=1]{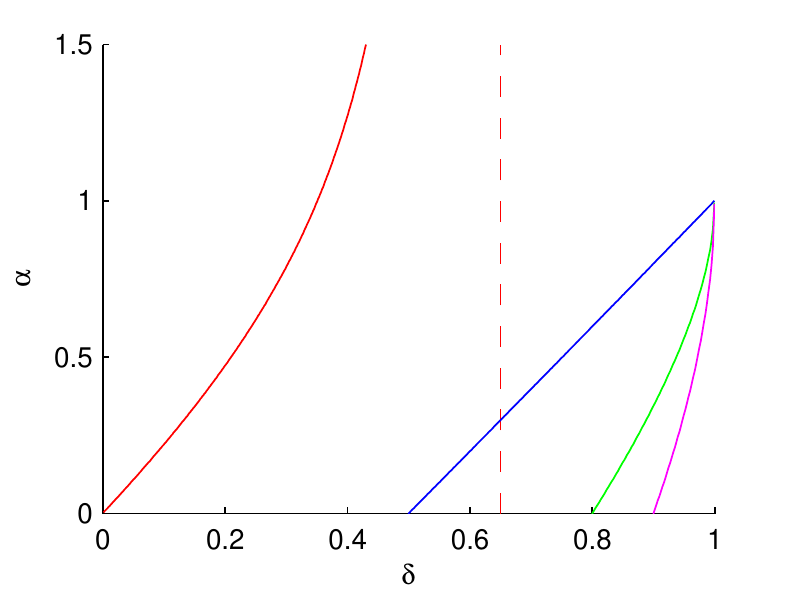}
                }
                \subfigure[$b=0.5$]    {
                        \label{Stability_curveb=05}
                       \includegraphics[scale=1]{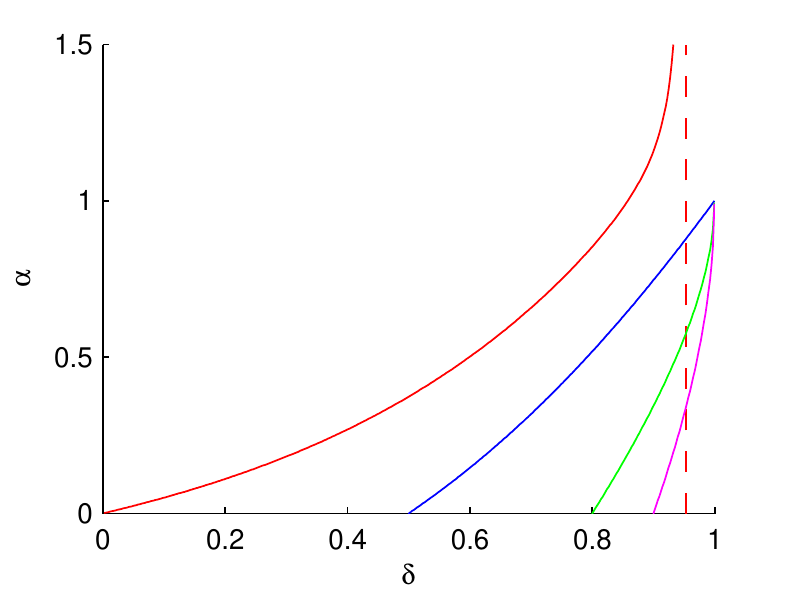}
               }
                \subfigure[$b=0.9$]    {
                        \label{Stability_curveb=09}
                        \includegraphics[scale=1]{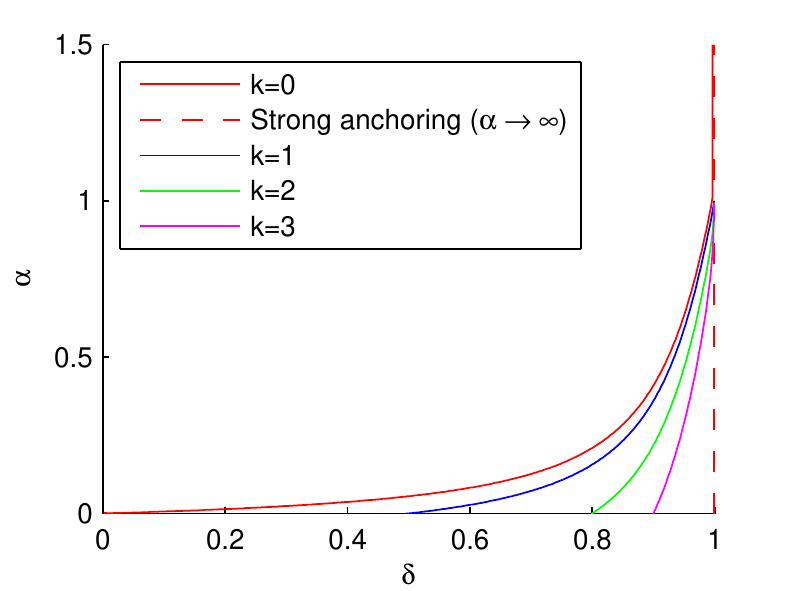}
                }
                \caption{The smallest solution to the compatibility condition (\ref{eq:s11}) and the strong anchoring limit, and the only solution to (\ref{eqn:s11.1}) for $k=1$, 2, 3 in the $(\delta, \alpha)$-plane for fixed $b$.}
                \label{fig:weak}
        \end{center}
\end{figure}

For $k>1$, the compatibility condition (\ref{eqn:s11.1}) has solutions if and only if
\begin{equation}
\frac{(1-\delta) \alpha (1+b)}{\delta \alpha + \alpha^2 b - \alpha b \delta - \delta + k^2 ( 1-\delta)} \sqrt {
\frac {{k}^{2}-\delta}{1-\delta}}  <0.
\end{equation} In particular, this implies
$\frac{k^2 + \alpha^2 b}{\alpha b - \alpha + 1 + k^2} <
\delta_{1,k}< 1
$, which in turn requires that $\alpha < 1$.
Therefore, the compatibility condition (\ref{eqn:s11.1}) only admits azimuthal perturbations with $k \geq 1$ if $0 < \alpha < 1$.
In Figure~\ref{fig:weak}, we plot the stability curves in the $(\delta,\alpha)$-plane, for fixed values of $k$ and $b$. These curves define the stability of the defect-free state with respect to the different kinds of perturbations e.g.\ the case $k=0$ corresponds to purely radial perturbations whereas $k\geq 1$ describe perturbations with azimuthal dependence. As $\delta$ increases, the system loses stability upon crossing one of the curves and as is evident from Figure~\ref{fig:weak}, larger values of $\alpha$ stabilize the defect-free state with respect to all perturbations. Similar comments apply to the confinement parameter $b$. As $b$ increases, the defect-free state is stable with respect to larger class of perturbations, as can be seen from Figure~\ref{fig:weak}. 
The critical value $\alpha = 1$ corresponds to $\xi= R_{\text{outer}}$ and the defect-free state is stable with respect to all azimuthal or symmetry-breaking perturbations for $\alpha>1$ or $\xi > R_{\text{outer}}$.

\section{Nematic Equilibria with Defects}
\label{sec:defects}
In this section, we consider nematic equilibria on a 2D annulus with regularly spaced boundary defects in the Oseen-Frank framework. 

To model states with regularly spaced defects along $\partial \Omega$, we split $\Omega$ into $N$ sectors with $N \in \mathbb{N}$, and define the region $\Omega_N$ to be
\begin{equation} \label{eq:d2}
\Omega_N = \left\lbrace (r,\phi) \in \mathbb{R}^2 : b \leq r \leq 1, 0 \leq \phi \leq \frac{2\pi}{N} \right\rbrace.
\end{equation}

\begin{figure}[!ht]
        \begin{center}
        \subfigure[Rotated state 1 ($U_1$)]    {
                        \label{Rotated_state_1}
                        \includegraphics[scale=0.2]{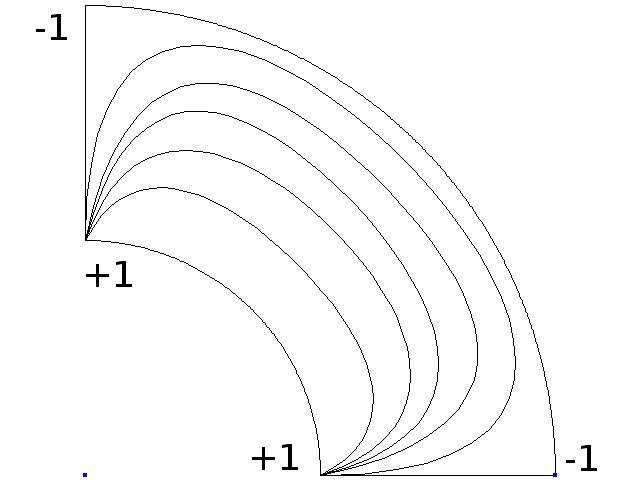}
                }
                \subfigure[Diagonal state 1 ($D$)]    {
                        \label{Diagonal_state_1}
                        \includegraphics[scale=0.2]{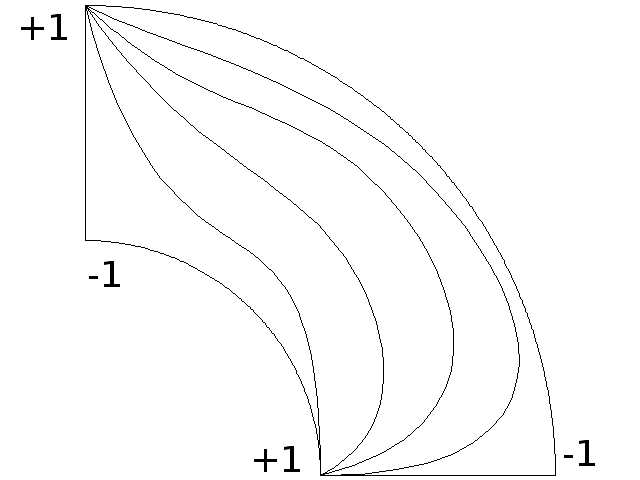}
               }
                \subfigure[Rotated state 2 ($U_2$)]    {
                        \label{Rotated_state_2}
                        \includegraphics[scale=0.2]{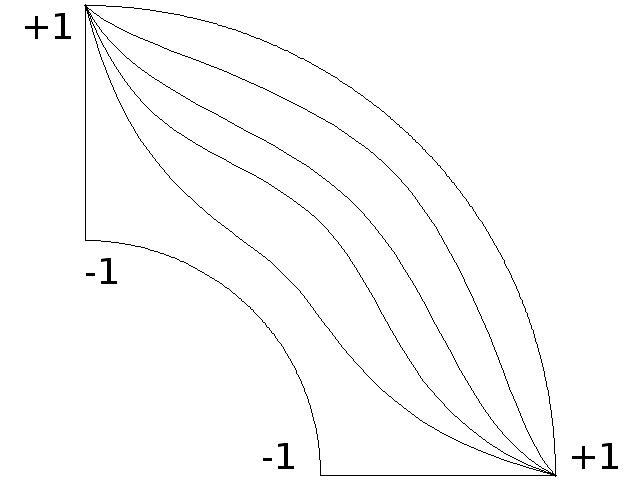}
                }
                \subfigure[Rotated state 3 ($U_3$)]    {
                        \label{Rotated_state_3}
                        \includegraphics[scale=0.2]{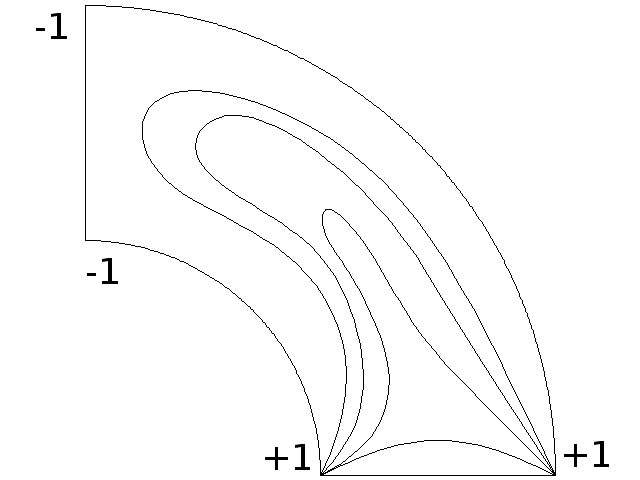}
                }
                \caption{The four distinct defect arrangements with the sector $\Omega_4$, with the strength of the defects denoted by the numbers in the corners, and sketches of the director.}
                
                \label{DefectConfigurations}
        \end{center}
\end{figure}

We assume strong tangent boundary conditions on all four edges of $\Omega_N$, $r=b$ and $r=1$, $\phi=0$ and $\phi=\frac{2 \pi}{N}$, so that there are necessarily discontinuities at the four corners. We only consider states with boundary defects of strength $m =\pm 1$, states with defects of higher strength are expected to have higher energy \cite{dg}.The $m=+1$ defect has a local radial splay profile and the $m=-1$ defect has a local bend profile.
There are 4 distinct arrangements of defects on the boundary, $\partial\Omega_N$, as shown in Figure \ref{DefectConfigurations}, to be considered. Three of these arrangements correspond to generalizations of the rotated state found in a square and rectangular wells, denoted $U_1$, $U_2$ and $U_3$, where the director field connects two adjacent $+1$-defects along an edge, and the fourth arrangement is a generalized diagonal state, $D$, for which the director field connects two diagonally opposite $+1$-defects. In the $U_1$ state, the two $+1$-defects lie along the edge $r=b$; the two $+1$-defects lie at the vertices of the edge $r=1$ for the rotated state $U_2$ and the two $+1$-defects are located at the vertices of the edge, $\phi=0$, for $U_3$. Other rotated and diagonal states are rotationally equivalent to one of the four cases enumerated in Figure \ref{DefectConfigurations}. 
The overall configuration in $\Omega$ is a superposition of the different states in the different sectors, $\Omega_N$. We note that for odd $N$, a superposition of $D$ or $U_3$ states would be physically unrealistic as it would require the superposition of a $+1$ and a $-1$-type defect at one of the vertices.

\subsection{The One Constant Approximation}
Using the one-constant approximation, for which $K_1 = K_3 \equiv K$, the OF energy reduces to
\begin{equation}
\label{eq:d1}
E[\theta]: = \iint_{\Omega} \frac{K}{2} \left|\grad \theta \right|^2 dA.
\end{equation}

The key step is to construct solutions of the Laplace equation, $ \grad^2 \theta = 0$ on $\Omega_N$, with Dirichlet boundary conditions (determined by the tangent conditions) on the four edges. The tangent conditions require that $\theta = \phi \pm \frac{\pi}{2}$ on $r=b$ and $r=1$; $\theta = 0$ on $ \phi=0$ and $\theta = \left\{\frac{2\pi}{N}, \frac{2\pi}{N} \pm \pi \right\}$ on the edge $\phi = \frac{2 \pi}{N}$.  Any admissible $\theta$, subject to these Dirichlet boundary conditions, can be written as
\begin{equation} \label{eq:d3}
\theta = a_0 \phi + a_1 f_1 + a_2 f_2 + a_3 f_3 + a_4 f_4,
\end{equation}
where the canonical functions $f_i$ are solutions of the Laplace equation with boundary conditions:
\begin{itemize}
\item $f_1 (r,0) = f_1 (r,2\pi/N) =0$, $f_1(b,\phi) = 0$ and $f_1(1,\phi) = 1$

\item $f_2 (r,0) = f_2 (r,2\pi/N) =0$, $f_2(b,\phi) = 0$ and $f_2(1,\phi) = \phi$

\item $f_3 (r,0) = f_3 (r,2\pi/N) =0$, $f_3(b,\phi) = 1$ and $f_3(1,\phi) = 0$

\item $f_4 (r,0) = f_4 (r,2\pi/N) =0$, $f_4(b,\phi) = \phi$ and $f_4(1,\phi) = 0$ .
\end{itemize}
We compute the canonical functions, $f_1 \ldots f_4$, using separation of variables, to be:
\begin{align} \label{eq:d4}
f_1 (r,\phi) &= \sum^{\infty}_{n=1} \frac{4 \sin \left(\frac{(2n-1) N}{2} \phi \right) \left( \frac{r^{-\frac{(2n-1)N}{2}} - b^{-(2n-1)N} r^{\frac{(2n-1)N}{2}}}{1 -  b^{-(2n-1)N}} \right)}{(2n-1)\pi},  \\
f_2 (r,\phi) &= \sum^{\infty}_{n=1} \frac{4 (-1)^{n+1} \sin \left(\frac{N n }{2} \phi \right) \left( \frac{r^{-\frac{n N}{2}} - b^{-nN} r^{\frac{nN}{2}}}{1-b^{-nN}} \right)}{N n}, \\
f_3 (r,\phi) &= \sum^{\infty}_{n=1} \frac{4 \sin \left(\frac{(2n-1) N}{2} \phi \right) \left( \frac{r^{-\frac{(2n-1)N}{2}} -  r^{\frac{(2n-1)N}{2}}}{b^{-\frac{(2n-1)N}{2}}- b^{\frac{(2n-1)N}{2}}} \right)}{(2n-1)\pi}, \\
f_4 (r,\phi) &= \sum^{\infty}_{n=1} \frac{4 (-1)^{n+1}\sin \left(\frac{N n \pi}{2} \phi \right) \left(\frac{r^{-\frac{n N}{2}} - r^{\frac{n N}{2}}}{b^{-\frac{n N}{2}} - b^{\frac{n N}{2}}} \right)}{N n} .
\end{align}
\begin{figure}[!ht]
        \begin{center}
        \subfigure[Rotated state 1 ($U_1$)]    {
                        \label{Rotated state 1}
                        \includegraphics[scale=0.5]{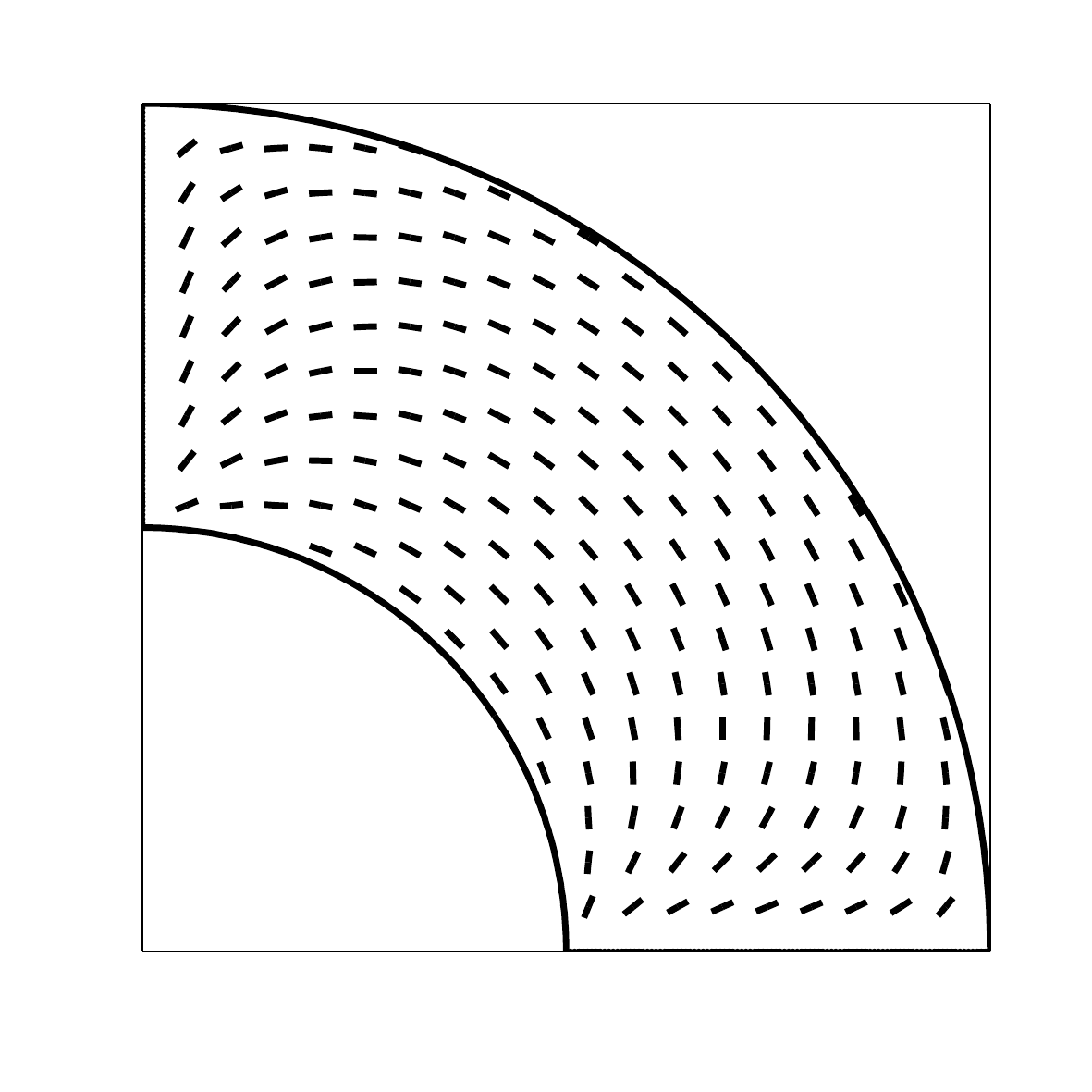}
                }
                \subfigure[Diagonal state 1 ($D$)]    {
                        \label{Diagonal state 1}
                        \includegraphics[scale=0.5]{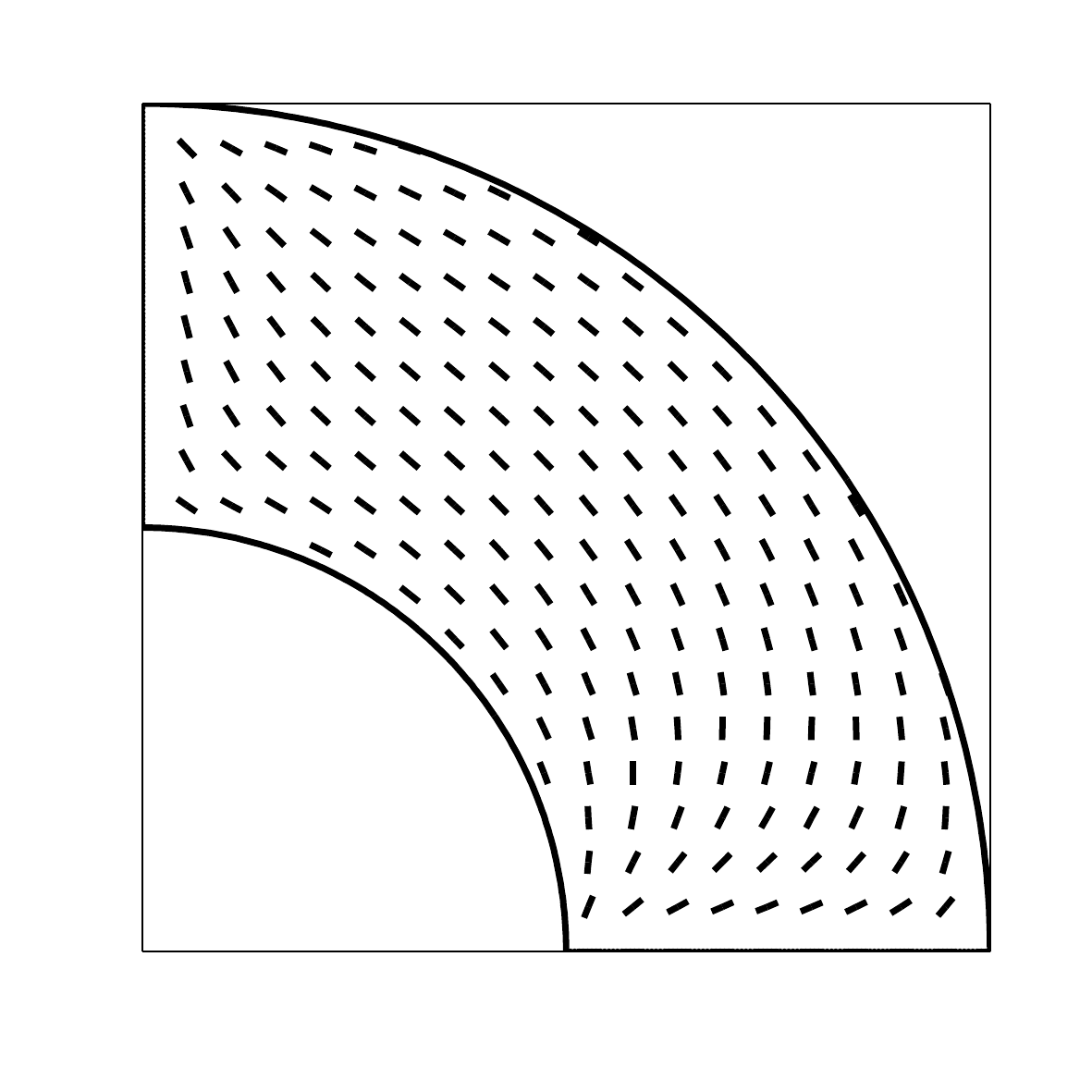}
               }
                \subfigure[Rotated state 2 ($U_2$)]    {
                        \label{Rotated state 2}
                        \includegraphics[scale=0.5]{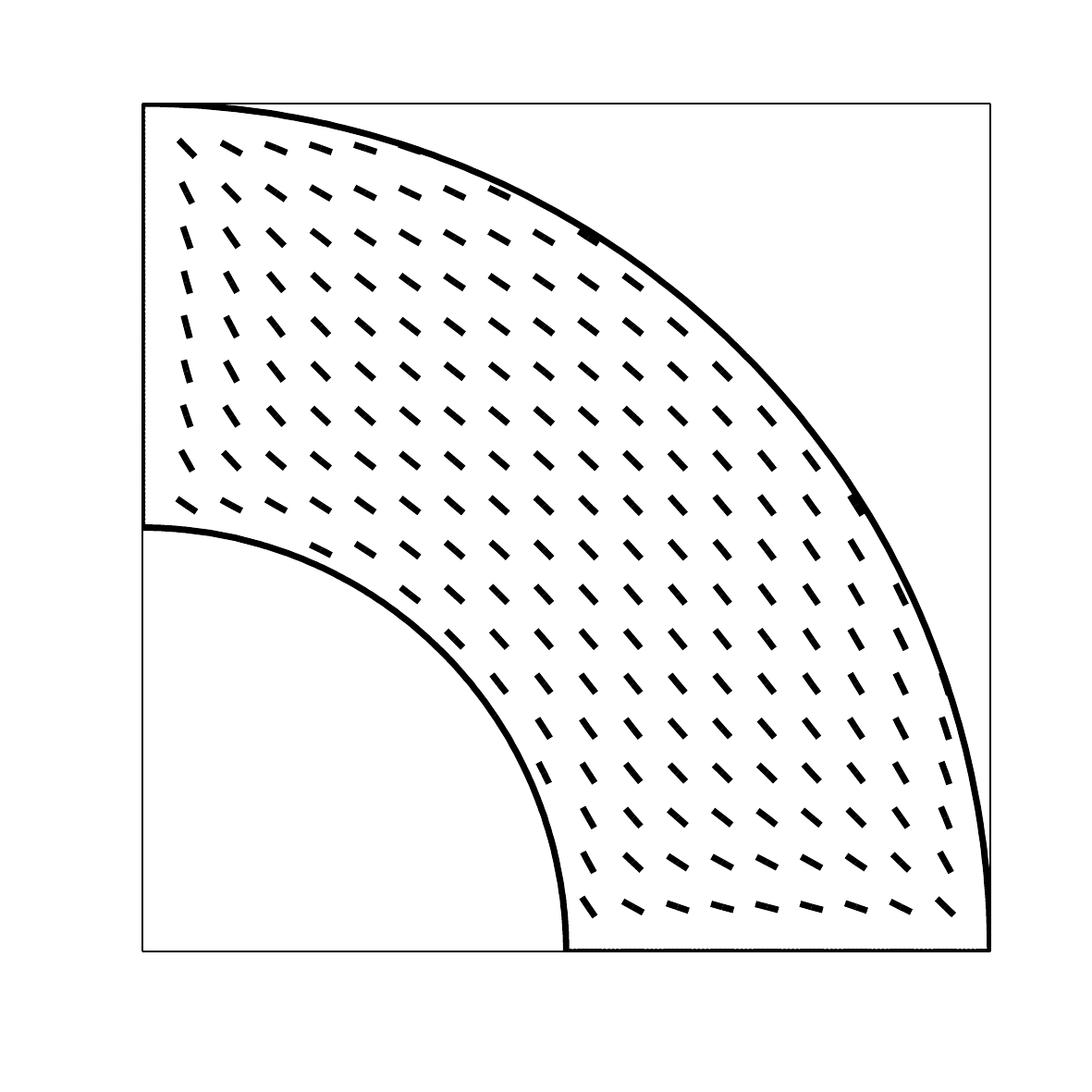}
                }
                \subfigure[Rotated state 3 ($U_3$)]    {
                        \label{Rotated state 3}
                        \includegraphics[scale=0.5]{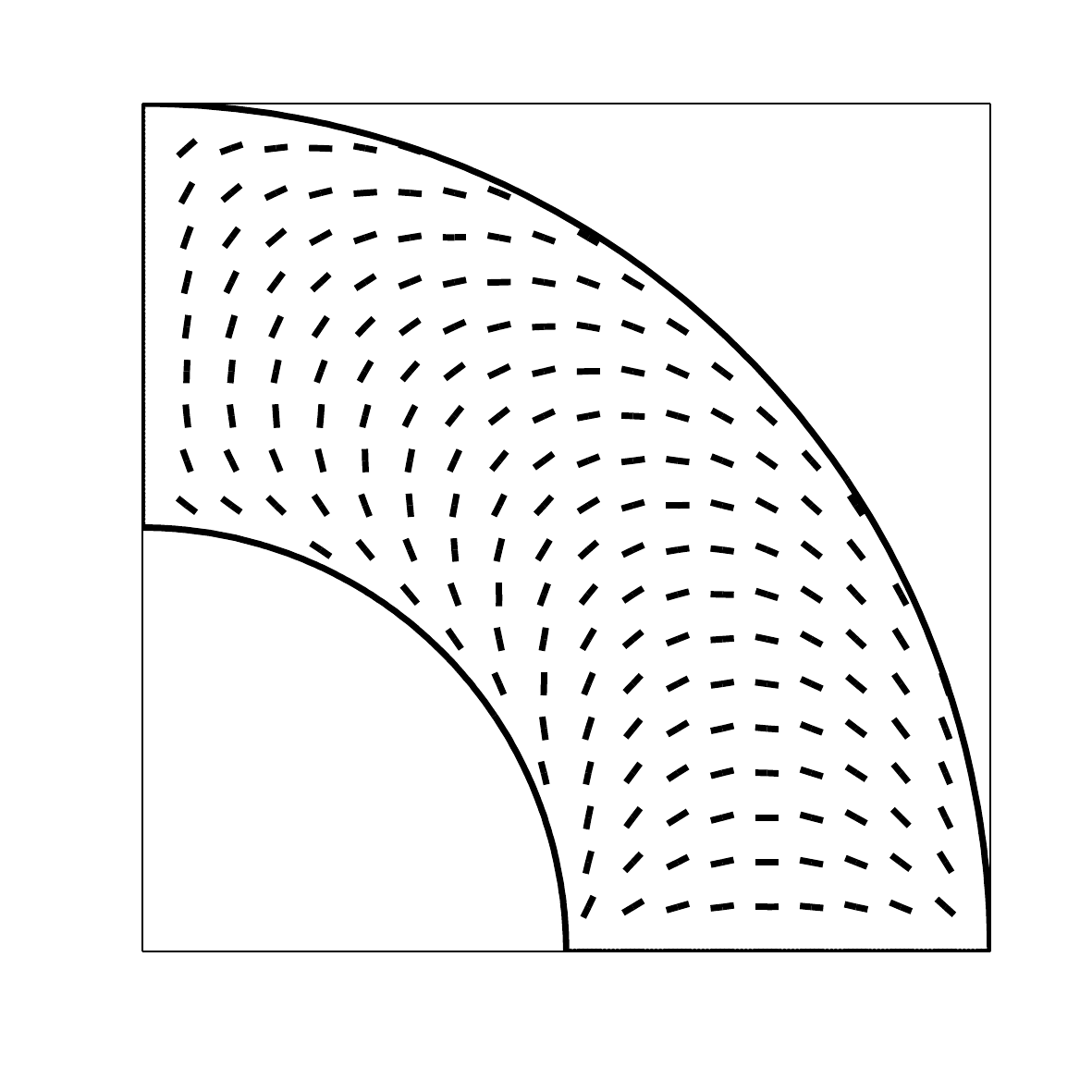}
                }
                \caption{The director of the four states within $\Omega_4$ and $b=0.5$, in the one constant approximation.}
                \label{AnnulusDirectorPlots}
        \end{center}
\end{figure}

Sample director plots using the canonical functions are shown in Figure \ref{AnnulusDirectorPlots} for $N = 4$ and $b = 0.5$. Due to the discontinuities in $\theta$ at the corners of $\Omega_N$, we regularize the domain by removing disks of radius $\epsilon$ about each of the defects (to prevent the energy (\ref{eq:d1}) from diverging) and the new regularized domain is denoted by $\Omega_{N_\epsilon}$. The length $\epsilon$ is proportional to the defect core size \cite{The_Static_and_Dynamic_Continuum_Theory} and we assume that $\epsilon \ll b$. The boundary, $\partial \Omega_{N_\epsilon}$, consists of four straight edges, $C_1 \ldots C_4$, and four curved arcs of radius $\epsilon$ enclosing the corners denoted by $\gamma_1 \ldots \gamma_4$ respectively.
We use Green's Theorem to write the energy as
\begin{equation} \label{eq:d5}
E = \frac{K}{2} \iint_{\Omega_{N\epsilon}} \vert \nabla \theta \vert^2 d A = \frac{K}{2} \oint_{\partial \Omega_{N_\epsilon}} \theta \nabla \theta \cdot \mathbf{v} d s = \sum^{4}_{i=1} \frac{K}{2} \int_{C_i} \theta \nabla \theta \cdot \mathbf{v} d s + \frac{K}{2} \int_{\gamma_i} \theta \nabla \theta \cdot \mathbf{v} d s,
\end{equation} where $\mathbf{v}$ is the outward pointing normal to each boundary segment.
The contributions from the curved arcs, $\gamma_i$, are all of order $\epsilon$.
We define four functions related to the integral contributions from the straight edges;
\begin{align}\label{eq:d9}
s_1 (N,b) &:= 8 \sum^{\infty}_{n=1} \frac{\coth \left( \frac{N}{2} (2n-1) \ln b \right) +1}{2n-1}, \\
s_2 (N,b) &:= 8 \sum^{\infty}_{n=1} \frac{\text{csch} \left( \frac{N}{2} (2n-1) \ln b \right) }{2n-1}, \\
s_3 (N,b) &:= 8 \sum^{\infty}_{n=1} \frac{\coth \left( \frac{N}{2} n \ln b \right) +1}{n} , \\
s_4 (N,b) &:= 8 \sum^{\infty}_{n=1} \frac{\text{csch} \left( \frac{N}{2} n \ln b \right) }{n} .
\end{align}
As $\eps \to 0$, the regularized energies within $\Omega_{N_\epsilon}$ are then given by
\begin{equation}
E \sim K \pi \left(\log \left(\frac{1}{\epsilon}\right) + \tilde{E}\right) + O(\epsilon) ,
\end{equation}
where the normalized energy, $\tilde{E}$, is the interior distortion energy and the logarithmic term originates from the defects.
The normalized energies of the four states are:
\begin{align}\label{eq:d10}
\tilde{E}_{U_1} =& \frac{ s_1 (N,b) + s_4(N,b) - s_2(N,b) - s_3(N,b) }{4} + \frac{(N+2)^2}{4N} \log\left(\frac{1}{b}\right) + \frac{1}{2} \log \left( \frac{b}{N^2} \right), \\
\tilde{E}_{U_2} =& \frac{ s_1 (N,b) + s_4(N,b) - s_2(N,b) - s_3(N,b) }{4} + \frac{(N-2)^2}{4N} \log\left(\frac{1}{b}\right) + \frac{1}{2} \log \left( \frac{b}{N^2} \right), \\
\tilde{E}_{U_3} =& -\frac{s_1(N,b) + s_2(N,b)}{4} + \frac{1}{N} \log\left(\frac{1}{b}\right) + \frac{1}{2} \log\left(\frac{16b}{N^2}\right) ,\\
\tilde{E}_{D} =& \frac{s_2(N,b) - s_1(N,b)}{4} + \frac{1}{N} \log\left(\frac{1}{b}\right) + \frac{1}{2} \log\left(\frac{16b}{N^2}\right).
\end{align}

We make some immediate comments. In the case of a square or rectangular domain, the curved arcs, $\gamma_i$, make an order $\epsilon^2$ contribution to the OF energy \cite{lewis1}. In the case of an annulus, the curved arcs $\gamma_i$ make an order $\epsilon$ energy contribution due to the curvature of the boundary.

\begin{figure}[!ht]
\center
\includegraphics[scale=0.75]{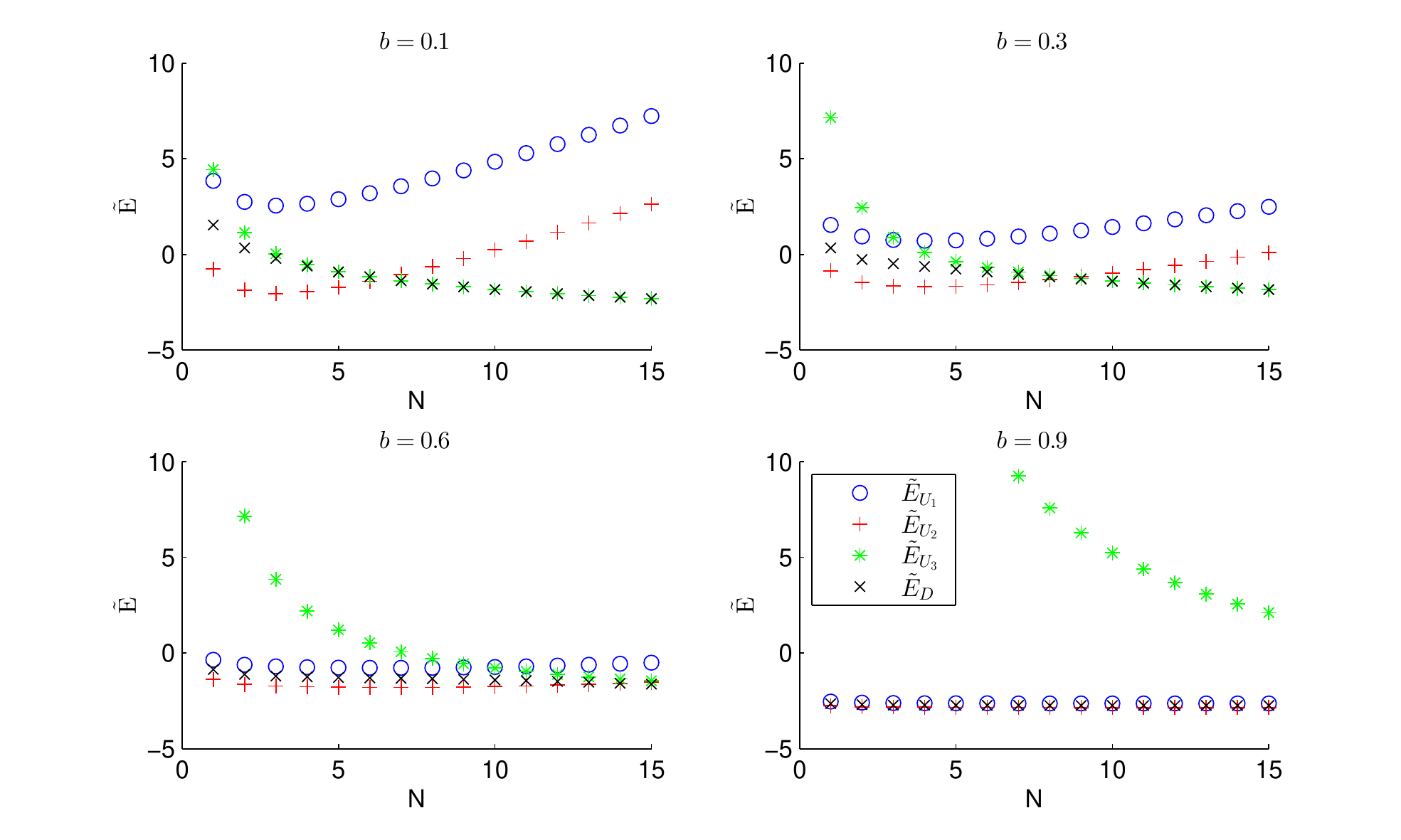}
\caption{The normalized energy of the four states plotted against $N$, for fixed values of $b$.}
\label{EnergyPlots}
\end{figure}

In Figure~\ref{EnergyPlots}, we plot the normalized OF energy of the four different states, as a function of $N$, in a fixed sector $\Omega_{N_\epsilon}$, for some fixed values of $b$. For small values of $N$, the $U_2$ state is the minimum energy state in the set $\left\{U_1, U_2, U_3, D \right\}$. In contrast, for a square or rectangle, the diagonal state has the minimum normalized energy \cite{lewis1}. As $N$ increases, there is a cross-over between the diagonal and rotated states and the $D$ state has the minimum normalized energy state for large $N$. This is consistent with the fact that $\Omega_N$ approaches a rectangle as $N \to \infty$. The cross-over critical value, $N = N_c$, depends on $b$; in particular, our limited simulations suggest that $N_c$ is an increasing function of $b$. This gives insight into preferred defect locations as a function of the annular aspect ratio.


\begin{figure}
\centering
\includegraphics[scale=1]{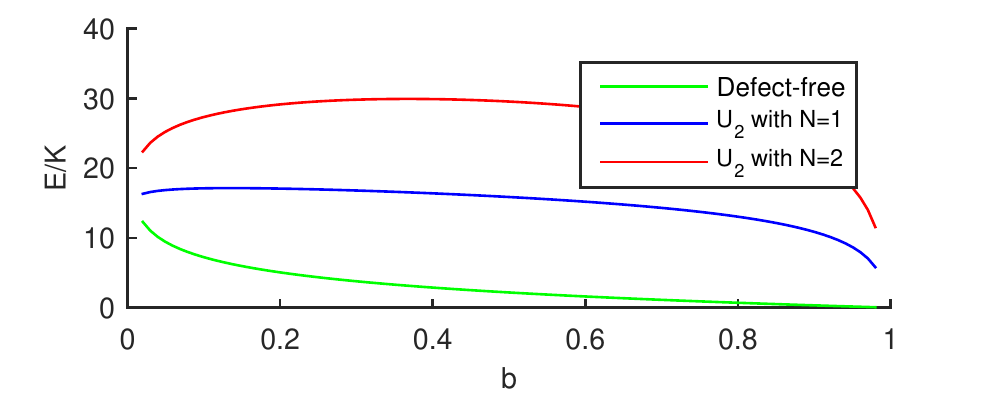}
\caption{The one constant Oseen-Frank energy of the three lowest energy states; the defect-free state and the $U_2$ states with $N=1,2$, for defect core size $\epsilon = 0.002$.}
\label{fig:IsotropicEnergy}
\end{figure}

In all cases, the defect-free state has lower energy than the generalized diagonal or rotated solutions, $\left\{U_1, U_2, U_3, D \right\}$, for small values of $\epsilon$, see Figure \ref{fig:IsotropicEnergy}. This stems from the $N\log(1/\epsilon)$- energy contribution from defects, which is absent in the defect-free case.

\subsection{Elastic Anisotropy}
\begin{figure}[!ht]
        \begin{center}
        \includegraphics[scale=0.75]{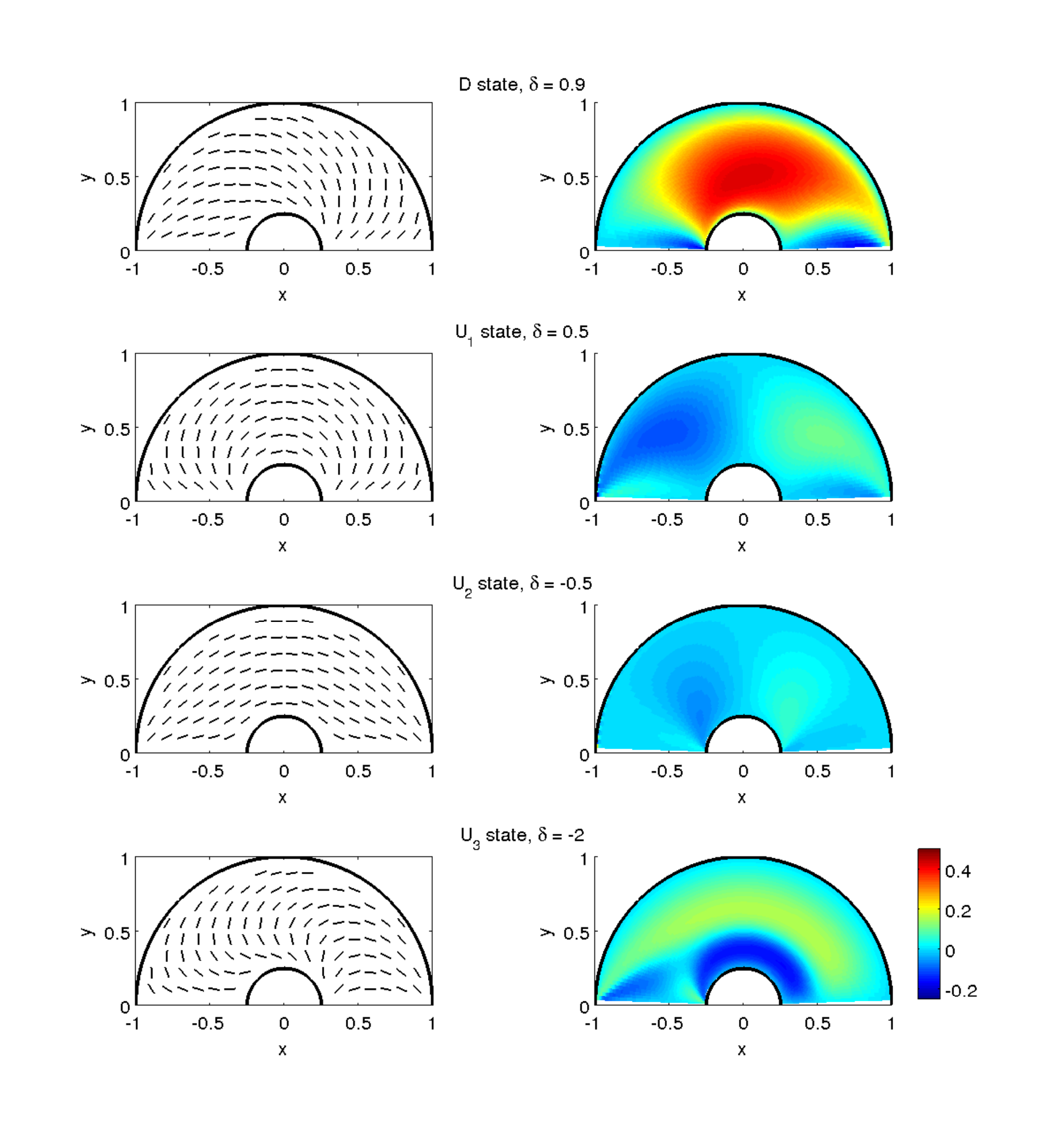}
\end{center}
\caption{Examples of the effect of elastic anisotropy on minimizers of the Oseen-Frank energy in $\Omega_2$ and $b=0.25$. We show numerically calculated director fields (left) and the difference in radians between the numerically evaluated $\theta$ and the one constant solution given in Equation (\ref{eq:d3}) (right).}
                \label{fig:AnisotropicDirector}
\end{figure}

We use the commercial PDE solver, COMSOL, to numerically solve the
Euler-Lagrange equation (\ref{eq:7}) with $\delta \neq 0$. We find
all $4$ states, $\left\{ U_1, U_2, U_3, D \right\}$, within a sector $\Omega_N$, see Figure
\ref{fig:AnisotropicDirector}. We further plot the structural
differences between the one-constant case (with $\delta=0$) and
the anisotropic cases in Figure \ref{fig:AnisotropicDirector}. The
structural details in the director are qualitatively the same.

\begin{figure}[!ht]
\centering
\includegraphics[scale=0.75]{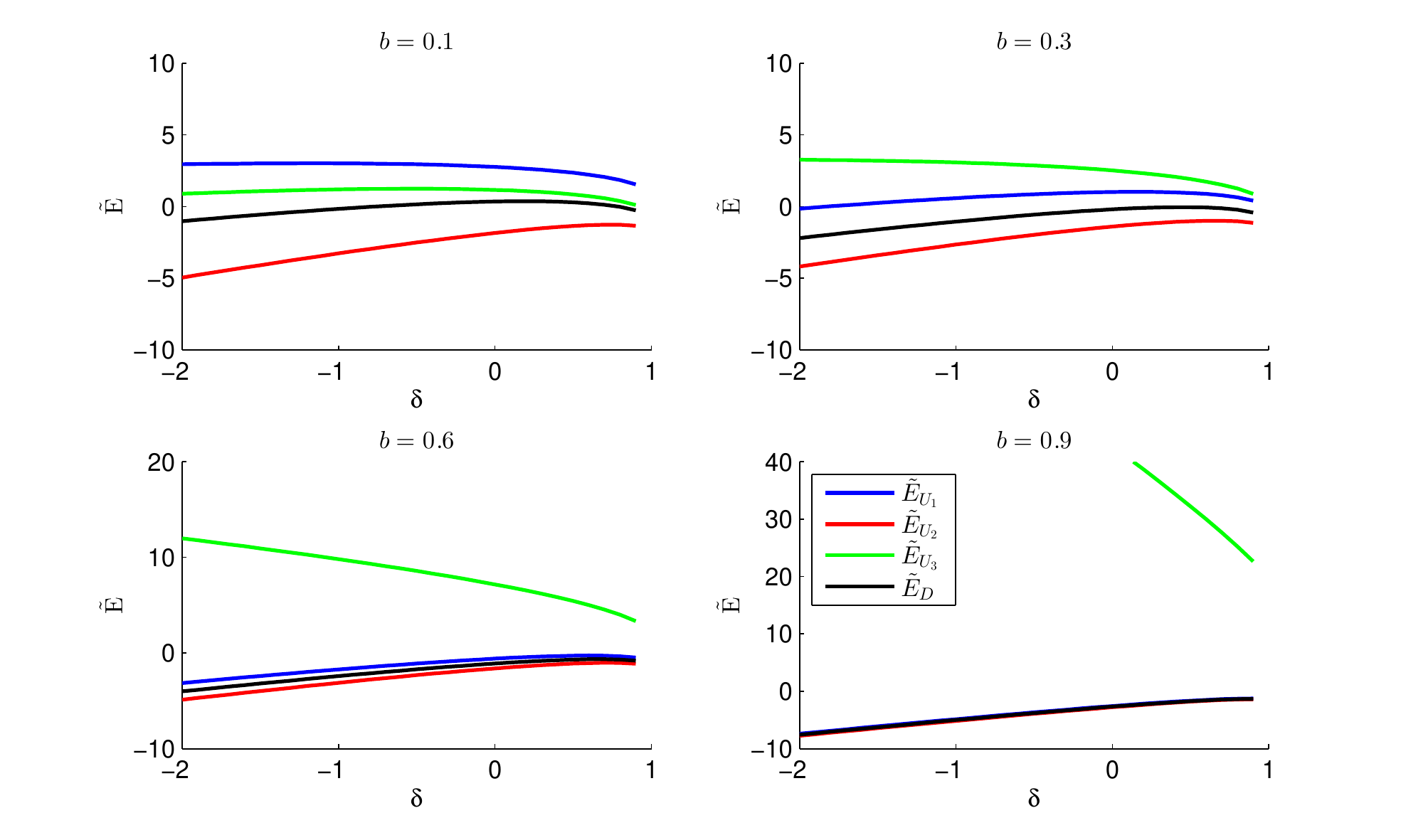}
              \caption{The normalized energies of the four states in $\Omega_2$, for fixed values of $b$.}
               \label{fig:AnisotropicEnergies}
\end{figure}

The regularized OF energy within $\Omega_N$ can be expressed as
the sum of the defect contributions and the normalized energy,
$\tilde{E}(\delta,b,N)$.
\begin{equation}
E \sim K_3 \pi \left( \left( 1 - \frac{3\delta}{4} \right)  \log \left(\frac{1}{\epsilon}\right) + \tilde{E}(\delta,b,N) \right) + O(\epsilon).
\end{equation}
In Figure \ref{fig:AnisotropicEnergies}, we plot the numerically evaluated normalized energies for fixed $b$ and $N$. For small $N$, the rotated $U_2$ state has the minimum
energy in the set, $\left\{U_1, U_2, U_3, D \right\}$, for all
admissible values of $\delta$ and $b$. As $b \to 1$, the
normalized energies of $U_1, U_2, D$ converge whereas the
normalized energy of the $U_3$ state diverges as $b \to 1$. From
Figure~\ref{AnnulusDirectorPlots}, the $U_3$-director is
constrained to rotate by approximately $\pi$ radians over a
distance of $(1-b)$ units, leading to the energy blow-up as $b \to
1$.

\begin{figure}
\centering
\includegraphics[scale=1]{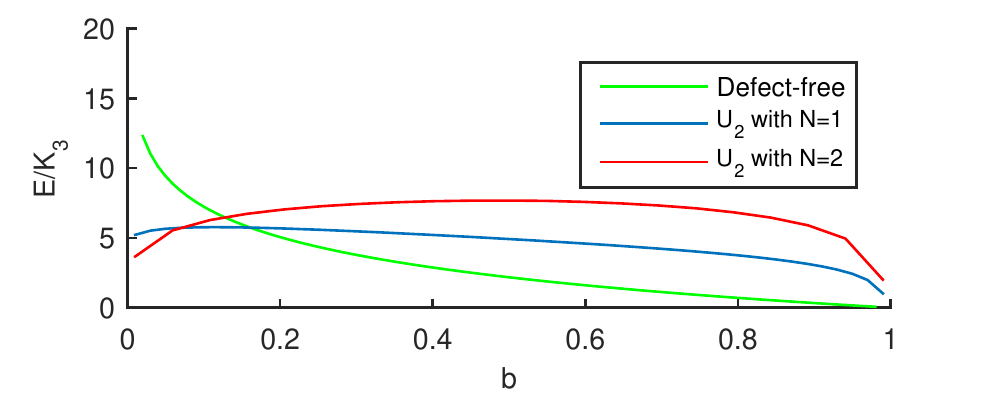}
\caption{The Oseen-Frank energy, with $\delta = 0.9$, of the three lowest energy states; the defect-free state and the $U_2$ states with $N=1,2$, for defect core size $\epsilon = 0.002$.}
\label{fig:AnisotropicEnergy}
\end{figure}

For suitable values of $\epsilon$ and $\delta$ close to 1, the defect-free state ceases to be the global minimizer of the Oseen-Frank energy, see Figure \ref{fig:AnisotropicEnergy}. As $b$ decreases, we see the $U_2$ state with $N=1$ minimizes the energy, than as $b$ decreases further the $U_2$ state with $N=2$ becomes the minimizing director field. 
This cross-over in energies demonstrates the importance of elastic anisotropy and has implications for experimental observations, \cite{lewis1}.




\section{The Defect-Free State in the Landau-de Gennes Model}
\label{sec:mironescu}

This section focuses on the defect-free state in the
two-dimensional LdG theory. We study maps $\Qvec:
\Omega \to S_2$ where $\Omega = \left\{ \xvec \in \Rr^2: b \leq
|\xvec| \leq 1 \right\}$ and $S_2$ is the space of symmetric
traceless $2\times 2$ matrices, subject to the fixed boundary
condition
\begin{equation}
\label{eq:Qb}
\Qvec_b = \frac{1}{\sqrt{2}}\left(\nvec\otimes \nvec - \mvec\otimes \mvec \right)
\end{equation} on $r=b$ and $r=1$.
In (\ref{eq:Qb}), $\nvec = \left(\cos \phi, \sin\phi\right)$,
$\mvec = \left(-\sin\phi, \cos\phi \right)$ and $\phi$ is the
standard  polar angle defined by $\tan\phi = \frac{y}{x}$. We
define the defect-free state to be
\begin{equation}
\label{eq:l3}
\Qvec^* = s(r) \left(\nvec\otimes \nvec - \mvec\otimes \mvec \right),
\end{equation} where $s:\Omega \to \Rr$ is the unknown scalar order parameter
subject to the Dirichlet conditions, $s(b) = s (1) = \frac{1}{\sqrt{2}}$.
The LdG energy of $\Qvec^*$ can be computed using (\ref{eq:3}) and (\ref{eq:3b}) to be
\begin{equation}\label{eq:l4}
I[\Qvec^*] = 2\pi \int^{1}_{b} \left( (s')^2 + \frac{4s^2}{r^2} + \frac{\vert A \vert}{4L} \left(2s^2 - 1 \right)^2 \right) r d r
\end{equation}
and we define the optimal order parameter $s$ to be a minimizer of
(\ref{eq:l4}) subject to the Dirichlet conditions on $r=b$ and
$r=1$. It is straightforward to verify (see \cite{EJM:8458158,
lamy}) that the minimizing $s$ is a classical non-negative
solution of the following second-order ordinary differential equation
\begin{equation}
\label{eq:l5}
s'' + \frac{s'}{r} - \frac{4s}{r^2} = \frac{\left| A \right|}{L} s (2s^2 -1).
\end{equation} subject to  $s(b) = s (1) = \frac{1}{\sqrt{2}}$. If $s$ were
 negative for $r_0 < r < r_1$, then one could define the admissible
 comparison map, $w = |s|$ which would have the same energy as $s$ i.e.~ $I[w] = I[s]$ but
 would have discontinuous derivatives at points where $s$ vanishes,
 contradicting the regularity properties of minimizing solutions.
 One can check that $\Qvec^*$, thus defined, is indeed a critical point of the
 LdG energy in (\ref{eq:3}) and the second variation of the
 LdG energy about $\mathbf{Q}^*$ is given by
\begin{equation}
\label{eq:l6}
\delta^2 I [\Qvec^*]:=\iint_{A} \frac{1}{2} \vert \nabla \Vvec \vert^2 + \frac{\vert A \vert}{2L} (2s^2 - 1) \vert \Vvec \vert^2  + \frac{\vert A \vert}{L} (\Qvec^* \cdot \Vvec)^2 d A,
\end{equation} where $\Vvec(r, \phi) \in S_2$  is an
admissible perturbation about $\Qvec^*$ with $\Vvec(b,\phi) = \Vvec(1,\phi) = 0$.

We use the following basis for the space of symmetric, traceless $2\times 2$ matrices,
\begin{equation} \label{eq:l1}
\mathbf{E} =  \left( \begin{array}{cc}
\cos(2\phi) & \sin(2\phi) \\
\sin(2\phi) & -\cos(2\phi)  \\
 \end{array} \right)
\end{equation}
and
\begin{equation}\label{eq:l2}
\mathbf{F} =  \left( \begin{array}{cc}
-\sin(2\phi) & \cos(2\phi) \\
\cos(2\phi) & \sin(2\phi)  \\
 \end{array} \right).
\end{equation} Then, we can write an arbitrary $\Vvec$ as
\begin{equation}\label{eq:l7}
\Vvec = v(r,\phi) \mathbf{E} + w(r,\phi) \mathbf{F},
\end{equation}
where $v(b,\phi) = v(1,\phi) = w(b,\phi) = w(1,\phi) =0$ and $v, w$ are
$2\pi$-periodic in $\phi$. We use a basis decomposition for the functions,
$v$ and $w$, as shown below:
\begin{eqnarray}\label{eq:l8}
v(r,\phi) = \sum_{n=0}^\infty a_{n} (r) \cos(n\phi) + b_{n}(r) \sin(n\phi), \\
w(r,\phi) = \sum_{n=0}^\infty c_{n} (r) \cos(n\phi) + d_{n}(r) \sin(n\phi),
\end{eqnarray}
where $a_n = b_n = c_n = d_n =0$ on $r =b,1$ for all $n$.
The second variation can, thus, be expressed as
\begin{equation}\label{eq:l9}
\delta^2 I [\Qvec^*] = 2\pi  L_0[a_0,a_0,c_0,c_0] + 2\pi \sum^{\infty}_{n=1} L_n[a_n,b_n,c_n,d_n] ,
\end{equation}
where the functional $L_n$ is defined to be
\begin{align}
L_n [a,b,c,d] := & \int_{b}^{1} \left( (a')^2 + (b')^2 +(c')^2 + (d')^2 + \left(\frac{n^2 + 4}{r^2}\right) \left( a^2 + b^2 + c^2 + d^2 \right) \right. \nonumber \\
& \left. \quad + \frac{8n}{r^2} \left( a d - c d \right) \right) r + \frac{\vert A \vert}{L} \left( a^2 + b^2 + c^2 + d^2 \right) (2s^2 - 1) r +  \frac{4\vert A \vert}{L}s^2 \left( a^2 + b^2\right)r d r .
\end{align}

\paragraph{Proposition 2:} The defect-free state, $\Qvec^*$,
defined in (\ref{eq:l3}), is a locally stable equilibrium of the
LdG energy in (\ref{eq:3}) with $w(\Qvec, \grad \Qvec) =
\frac{L}{2}|\grad \Qvec|^2$, for all $0<b<1$ and $\frac{|A|}{L}>\frac{3 (b^2  +1)^2}{2 b^4}$.

\begin{proof} We follow the same strategy
as in \cite{mironescu}, where the author establishes the local
stability of the defect-free state on a disc, as opposed to an
annulus in two dimensions. For a disc, the order parameter is
defined for $0\leq r \leq 1$, the order parameter vanishes at $r=0$
and is a monotonically increasing function subject to a Dirichlet
condition at $r=1$. The function $s$, as defined in (\ref{eq:l5}),
cannot be a monotone function for $b\leq r \leq 1$ by virtue of
the imposed boundary conditions,  $s^2 \leq \frac{1}{2}$ by an
immediate application of the maximum principle and hence $s$ must
have an intermediate local minimum \cite{ejam2010majumdar}. The
key ingredient is to show that $\min L_n[a_n, b_n, c_n, d_n]>0$
for non-trivial $a_n, b_n, c_n, d_n$ and $\min L_0[a_0,a_0, c_0,c_0]>0$
for non-trivial $a_0$ and $c_0$.

\paragraph{Lemma 1:}
The second variation, $\delta^2 I[\Qvec^*]>0$ if $\min L_0>0$, $\min L_1 > 0$ and $\min L_2 >0$.
\begin{proof}
It suffices to show that for $n \geq 1$,
\begin{align}\label{eq:l10}
L_{n+2}(a,b,c,d) - L_n(a,b,c,d) &=
\int^{1}_{b} \left( \frac{4n+4}{r^2} (a^2 + b^2 + c^2 +d^2) + \frac{16}{r^2} (ad - bc) \right) r d r \geq 0
\end{align}
Using Young's inequality,
\begin{align}
(4n+4) (a^2 + b^2 + c^2 +d^2) + 16(ad - bc) 
& \geq (4n -4 ) (a^2 + b^2 + c^2 +d^2) \geq 0
\end{align}
\end{proof}
The problem of stability now reduces to a study of the functionals
$L_0$, $L_1$ and $L_2$.
The integral $L_0$ is given by
\begin{equation}
\label{eq:L0}
L_0: = \int_{b}^{1} \left[ \left(\frac{d a_0}{dr}\right)^2 + \left(\frac{d c_0}{dr} \right)^2 \right]  + \frac{4}{r}\left(a_0^2 + c_0^2 \right) r~ dr +  \int_{b}^{1} \frac{|A|}{L}(2s^2 - 1) (a_0^2 + c_0^2)r + \frac{4 |A|}{L}s^2 a_0^2 r~ dr.
\end{equation}
It suffices to show that
$L_0\geq I_0 = \int_{b}^{1} \left(\frac{d a_0}{dr}\right)^2 r + \frac{4}{r}a_0^2 + \frac{|A|}{L}(2s^2 - 1) a_0^2 r dr >0$ for all admissible $a_0$ such that $a_0(b)=a_0(1) =
0.$ We define $a(r) := s(r) a_0(r)$ and recall the governing ordinary differential equation
for $s$ in (\ref{eq:l5}). Then
\begin{eqnarray}
\label{eq:I0}
I_0 = \int_{b}^{1} \frac{d}{dr}\left( s \frac{ds}{dr} a^2 r \right) + s^2 \left(\frac{da}{dr}\right)^2 r dr  = \int_{b}^{1}s^2 \left(\frac{da}{dr}\right)^2 r dr > 0
\end{eqnarray}
for any non-trivial $a: \left[b, 1 \right] \to \Rr$ with $a(b) = a(1) = 0$.

It is straightforward to check that
\begin{align}
\label{eq:l11}
 L_1\left(a_1, b_1, c_1, d_1 \right) \geq  L_1\left(\sqrt{a_1^2 + b_1^2},0, \sqrt{c_1^2 + d_1^2},0 \right)
 \end{align} by using the change of variable, $a_1 = \sqrt{a_1^2 +
 b_1^2}\cos \theta$, $b_1 = \sqrt{a_1^2 +
 b_1^2}\sin \theta$ in the definition of $L_1$ above. Let
 $a=\sqrt{a_1^2 +
 b_1^2}$ and $d = \sqrt{c_1^2 + d_1^2}$ where  $a(b) = a(1) = d(b) = d(1) = 0$ .
At this step, we diverge from \cite{mironescu}. We minimize the
functional (\ref{eq:l4}), subject to the boundary conditions
$s(b)=0$ and $s(1) = \frac{1}{\sqrt{2}}$ and define $u$ to be a
corresponding energy minimizer. The existence of such a function
$u$ is guaranteed from the direct methods in the calculus of
variations and the minimizer $u$ is a classical solution of
\begin{equation}\label{eq:l12}
u'' + \frac{u'}{r} - \frac{4u}{r^2} = \frac{\vert A \vert}{L} (2u^2 - 1)u
\end{equation}
with $u(b) =0$ and $u(1) = \frac{1}{\sqrt{2}}$.
\paragraph{Proposition 3:} The function $u$ is
unique, monotonically increasing and $u \leq s$ for all $r$.
\begin{proof} One can show that $u$ is unique and monotonically increasing
by an immediate adaptation of the arguments in \cite{lamy}, which are omitted here for brevity.
We show that $u \leq s$.
We assume for a contradiction that  $s(r) < u(r)$ for some $r \in (b,1)$.
Then, there must exist $r_1$ and $r_2$, where $b < r_1 < r_2 \leq 1$,
for which $s(r_1) - u(r_1) = 0$, $s'(r_1) - u'(r_1) < 0$,
$s(r_2) - u(r_2) = 0$, $s'(r_2) - u'(r_2) > 0$ and $s(r) - u(r) < 0$ for $r \in [r_1,r_2]$.
We multiply the differential equations for $s$ and $u$ by $ru$ and $rs$,
subtract and integrate over $(r_1,r_2)$ to find
\begin{equation}
\int^{r_2}_{r_1} ru s'' - rs u'' + us' - su' d r =
\int^{r_2}_{r_1} \frac{\vert A \vert}{L} us (2s^2 - 2u^2) r d r.
\end{equation}
The RHS is negative since $s - u < 0$. The LHS is
\begin{eqnarray}\label{eq:l13}
\text{LHS} &=& \int^{r_2}_{r_1} u (r s')' - s(r u')' d r =
(r u s' - r s u') \bigg \vert^{r_2}_{r_1} \nonumber
 \\
&=& r_2 s(r_2) (s'(r_2) - u'(r_2) ) - r_1 s(r_1) (s'(r_1) -
u'(r_2)) > 0,
\end{eqnarray} yielding the desired contradiction.
\end{proof}
Since $s \geq u$, we immediately have
\begin{align}\label{eq:l14}
 L_1(a,0,d,0) & \geq    \int^1_b  \left( (a')^2 + (d')^2 + \frac{5}{r^2} \left(a^2  + d^2 \right) - \frac{8}{r^2}  ad   \right) r \nonumber \\  & \qquad + \frac{\vert A \vert}{L} \left( a^2 + d^2 \right) (2u^2 - 1)r + \frac{4\vert A \vert}{L} a^2 u^2 r d r \equiv  m .
\end{align}
It remains to show that $ \min m>0$. The rest follows by analogy with \cite{mironescu} with some technical differences. We compute the Euler-Lagrange equations associated with $a$ and $d$:
\begin{eqnarray}\label{eq:l15}
a'' + \frac{a'}{r} - \frac{5a}{r^2} + \frac{4d}{r^2} = \frac{\vert A \vert}{L} a (6u^2 - 1) - a\min m, \nonumber \\
d'' + \frac{d'}{r} - \frac{5d}{r^2} + \frac{4a}{r^2} = \frac{\vert A \vert}{L} d (2u^2 - 1) - d\min m , \nonumber
\end{eqnarray}
where $\min m$ acts as a Lagrange multiplier. This system is
satisfied weakly, with $\min m = 0$ and no boundary conditions,
by the functions $A = u'$ and $D = 2u/r$. Multiplying the
Euler-Lagrange equations (\ref{eq:l15}) by $rA$ and $rD$
respectively and the system for $A$ and $D$ by $ra$ and $rd$
respectively followed by integrating over $(b, 1)$ yields

\begin{equation}
- \min m \int^1_b r(aA + dD) d r = \left( r (a' A + d' D) \right) \bigg \vert^1_b .
\end{equation}
The RHS is negative since $a,d \geq 0$ , $a'(1), d'(1) < 0$ and
$a'(b), d'(b) > 0$. Therefore, $\min m \geq 0$. If $\min m = 0$,
then we must have $(a'(1) A(1) +\sqrt{2} d'(1) ) = b a'(b) A(b)$.
This implies $d'(1) =0$ and $a'(b) u'(b) = a'(1) u'(1) =0$. For
$a'(1) u'(1)  = 0$ to hold, we must either have $a'(1) = 0$ (and
hence $a = d = 0$) or $u'(1) = 0$. From the ODE for $u$ in
(\ref{eq:l12}), we see that $u''(1) = 2\sqrt{2}, $ which implies
that $u(1)$ is a local minimum contradicting the fact that $u$ is
a monotonically increasing function for $b\leq r \leq 1$.

The last step is to show that $L_2$ is positive for large values of $\frac{\vert A \vert}{L}$. We first show that $L_2 > 0$ in the limit as $\frac{\vert A \vert}{L} \rightarrow \infty$ by following the strategy laid out in \cite{3dshells}, and then derive an explicit condition in terms of $b$ and $\frac{\vert A \vert}{L}$.
We note that
\begin{equation}
L_2 \geq  \int^1_b \left( a_2'^2 + d_2'^2 + \frac{8 (a_2 - d_2)^2}{r^2} + \frac{\vert A \vert}{L} \left( a_2^2 (6s^2 - 1) + d_2^2 (2s^2 - 1) \right) \right) r ~d r \equiv m_1 ,
\end{equation}
where $a_2$ and $d_2$ are positive functions of $r$.

The first stage is to show that for sufficiently large $\frac{\vert A \vert}{L}$, $\min s(r)$ is controlled by $\frac{\vert A \vert}{L}$.
\paragraph{Proposition 4:}
$s \rightarrow \frac{1}{\sqrt{2}}$ uniformly as $\frac{\vert A \vert}{L} \rightarrow \infty$.
\begin{proof}
The function $s$ attains its minimum at $r_* \in (b,1)$. We define $s(r_*) = s_{\min}$. Therefore, $s'(r_*) = 0$, $s''(r_*) \geq 0$ and
\begin{equation}
\frac{\vert A \vert}{L} s_{\min}(2s_{\min}^2 - 1) + \frac{4s_{\min}}{r_*^2} \geq 0
\end{equation}
Rearranging, 
\begin{equation}
s(r) \geq s_{\min} \geq \sqrt{\frac{1}{2} - \frac{2L}{\vert A \vert r_*^2}} \rightarrow \frac{1}{\sqrt{2}} \qquad \text{as $\frac{\vert A \vert}{L} \rightarrow \infty$}
\end{equation}
\end{proof}

We let $a_2(r) = s(r) a(r)$, $d_2(r) = s(r) d(r)$, and use Hardy's trick to write $m_1$ as
\begin{equation}
m_1 = \int^1_b s^2 \left( a'^2 + d'^2 + \frac{4 (a^2 + d^2)}{r^2} - \frac{16 ad}{r^2} + \frac{4 \vert A \vert}{L} s^2  a^2 \right) r ~d r
\end{equation}
and consider the terms $s^2 \left( \frac{4 (a^2 + d^2)}{r^2} - \frac{16 ad}{r^2}\right) + \frac{4 \vert A \vert}{L} s^4  a^2$. Then
\begin{equation}
s^2 \left( \frac{4 (a^2 + d^2)}{r^2} - \frac{16 ad}{r^2}\right) + \frac{4 \vert A \vert}{L} s_{\min}^4  a^2 \geq 
 \frac{4s_{\min}^2 (a^2 +d^2)}{r^2} - \frac{256a^2}{r^2} - \frac{d^2}{16r^2} + \frac{4 \vert A \vert}{L} s_{\min}^4 a^2.
\end{equation}
The coefficients of $d^2$ are positive if the order parameter $s_{\min} > \frac{1}{8}$ for all $b \in (0,1)$.
From Proposition 4, we can always find some $T_1$ such that $s_{\min} > \frac{1}{8}$ for all $\frac{\vert A \vert}{L} > T_1$.


The coefficients of $a^2$ are positive if $\frac{\vert A \vert}{L} > T_2$, where we define $T_2 \equiv \frac{64}{s_{\min}^2 b^2}$.
Then $L_2 > 0$ if $\frac{\vert A \vert}{L} > \max \left\lbrace T_1, T_2 \right\rbrace$.

Finally, we show that $L_2$ is positive if $\frac{\vert A \vert}{L}$ is greater than an explicit positive quantity dependent on $b$.
We let $a(r) = \alpha(r) \cos(\theta(r))$ and $d(r) = \alpha(r) \sin(\theta(r))$; then
\begin{equation}
m_1 \geq  \int^1_b s^2 \left( \alpha'^2  + \frac{4\alpha^2}{r^2} - \frac{8\alpha^2 \sin(2\theta)}{r^2} + \frac{4\vert A \vert}{L} s^2 \alpha^2 \sin^2(\theta) \right) r ~d r \equiv m_2 \label{eq:lowtempL2}.
\end{equation}
The minimum of $m_2$, with respect to $\theta$, occurs at
\begin{equation}
\theta = \frac{1}{2} \arctan\left( \frac{4L}{\vert A \vert r^2 s(r)^2} \right)
\end{equation}
and
\begin{equation}
m_2 = \int^1_b s^2 \left( \alpha'^2 + \frac{4\alpha^2}{r^2} - \frac{32 \vert A \vert \alpha^2}{r^4 L s^2 \sqrt{1 + \frac{16L^2}{\vert A \vert^2 s^4 r^4} }}  + 2s^2 \alpha^2 \frac{\vert A \vert}{L} \left( 1 - \frac{1}{\sqrt{1 + \frac{16L^2}{\vert A \vert^2 s^4 r^4} }} \right) \right) r ~d r .
\end{equation}
We note that $\alpha = 0$ is always a critical point of
\begin{equation}
\frac{4\alpha^2}{r^2} - \frac{32 \vert A \vert \alpha^2}{r^4 L s^2 \sqrt{1 + \frac{16L^2}{\vert A \vert^2 s^4 r^4} }}  + 2s^2 \alpha^2 \frac{\vert A \vert}{L} \left( 1 - \frac{1}{\sqrt{1 + \frac{16L^2}{\vert A \vert^2 s^4 r^4} }} \right)
\end{equation}
and is a minima if $\frac{\vert A \vert}{L} > \frac{3}{r^2 s^2}$. Using the fact that $s \geq \frac{\sqrt{2} b}{b^2 + 1}$ (see \cite{Golovaty}), then $L_2 > 0$ if
\begin{equation}
\frac{\vert A \vert}{L} >  \frac{3(b^2+ 1)^2}{2b^4} .
\end{equation}
\end{proof}

\section{Conclusions}
\label{sec:conclusion}
In this paper, we investigate nematic equilibria confined to two-dimensional annuli with strong or weak tangential anchoring modelled in the continuum Oseen-Frank and Landau-de Gennes theories. 

We compute explicit stability criteria for the ddefect-free state with a radially-invariant director field in both the Oseen-Frank and Landau-de Gennes theories. In the Oseen-Frank theory with strong anchoring, we demonstrate that the defect-free state undergoes a supercritical pitchfork bifurcation and compute a new energetically preferable defect-free state with spiral-like director profile. Furthermore, we demonstrate that this spiral-like state is locally stable for $\delta \rightarrow 1$. 
In the Oseen-Frank theory with weak anchoring, we produce a new stability diagram (Figure \ref{fig:weak}) in terms of the parameters $\delta$ and $\alpha$.
We demonstrate that the defect-free state is stable with respect to symmetry-breaking perturbations when the outer radius exceeds a critical value, which is precisely the surface extrapolation length of the nematic material in use. In parallel, we obtain an explicit local stability criterion for the defect-free state, in terms of the temperature and geometry, in the Landau-de Gennes framework.


We model nematic equilibria with defects in the Oseen-Frank theory on the boundary by computing local solutions of the Laplace equation on an annular sector with Dirichlet tangent boundary conditions. We find three rotated and one diagonal solution, by analogy with similar work done for squares and rectangles in \cite{mottramdavidson,lewis1,Tsakonas_Davidson_Brown_Mottram}.
We compute analytic expressions for the corresponding director fields and their one-constant Oseen-Frank energy and these energy expressions allow us to distinguish between the energy contributions from the defective corners and the bulk distortion and give information about the optimal number and arrangement of boundary defects as a function of geometry.
In particular, we find that the rotated $U_2$ state, which has been found experimentally in \cite{Oli_thesis, Jose_thesis}, has the minimum energy in the restricted class $\left\lbrace U_1,U_2, U_3, D \right\rbrace$, for small numbers of defects. 
In contrast, the diagonal solutions always have minimum energy for a square or a rectangle in the one-constant framework, see \cite{lewis1}. 
Furthermore, we investigate the effect of elastic anisotropy and demonstrate that for $K_3 \gg K_1$, the $U_2$ solution with boundary defects can have lower energy than the defect-free state for specific choices of the geometry. The singular limit $\delta \rightarrow 1$ will be investigated in future work.

\section{Acknowledgements}
We thank Oliver Dammone for valuable discussions. A.L.~ is supported by the Engineering and Physical Sciences Research Council (EPSRC) studentship. A.M.~ is supported by an EPSRC Career Acceleration Fellowship EP/J001686/1 and EP/J001686/2, an OCCAM Visiting Fellowship and the Keble Advanced Studies Centre.

\addcontentsline{toc}{chapter}{Bibliography}
\bibliographystyle{plain}
\bibliography{Literature_Review}

\end{document}